\documentclass{article}
\usepackage{amsmath}
\usepackage{amsfonts}
\usepackage{graphicx}
\usepackage[caption=false]{subfig}
\usepackage{cleveref}
\usepackage{physics}
\usepackage{color}

\newcommand*{\affaddr}[1]{#1} 
\newcommand*{\affmark}[1][*]{\textsuperscript{#1}}
\newcommand*{\email}[1]{\texttt{#1}}

\begin{document}

\title{Solving Differential Equation with Constrained Multilayer Feedforward
Network}

\author{%
    Zeyu Liu\affmark[1], Yantao Yang\affmark[1,2], and Qing-Dong
    Cai\affmark[1,2,3] \\
    \affaddr{\affmark[1]College of Engineering} \\
    \affaddr{\affmark[2]State Key Laboratory for Turbulence and Complex
    System}\\
    \affaddr{\affmark[3]Center for Applied Physics and Technology} \\
    \email{\{zeyuliu,yantao.yang,caiqd\}@pku.edu.cn} \\
    \affaddr{Peking University}
}

\date{\today}
\maketitle

\begin{abstract}
    In this paper, we present a novel framework to solve differential equations
    based on multilayer feedforward network. Previous works indicate that 
    solvers based on neural network have low accuracy due to that the boundary 
    conditions are not satisfied accurately. The boundary condition is now 
    inserted directly into the model as boundary term, and the model is a 
    combination of a boundary term and a multilayer feedforward network with 
    its weight function. As the boundary condition becomes predefined 
    constraintion in the model itself, the neural network is trained as an 
    unconstrained optimization problem. This leads to both ease of training and
    high accuracy. Due to universal convergency of multilayer feedforward 
    networks, the new method is a general approach in solving different types 
    of differential equations. Numerical examples solving ODEs and PDEs with 
    Dirichlet boundary condition are presented and discussed.
\end{abstract}

\section{Introduction}

Differential equations, including ordinary differential equations (ODEs) and 
partial differential equations (PDEs), are key mathematical models for various 
physics and engineering applications. In most situations, it is impractical to 
find analytical solutions, and numerical solutions become increasingly popular 
for these problems. When solving ODE/PDEs, one seeks for a function satisfying 
both (1) the differential equations within the domain, and (2) all 
initial/boundary conditions. Common numerical methods for ODEs are Runge-Kutta 
methods, linear multistep methods, and predictor-corrector methods 
\cite{burden2001numerical}. As for PDEs, tremendous methods for discretizing
the physical space or spectral space are developed, and the most common choices
are finite difference method (FDM), finite volume method (FVM), finite element
method (FEM), and spectral method. These methods are special cases of weighted 
residual method. Galerkin method is another numerical method based on weighted 
residual method for converting a continuous operator to a discrete form. It 
applies the method of variation of parameters to a function space and converts 
the original equation to a weak formulation. 

In the present study, we take advantage of the fast developing machine learning 
technique and propose a framework of solving ODE/PDE by applying the variation 
of parameters of a neural network \cite{raissi2017physics,lagaris1998artificial}.
Neural network (NN) is inspired by the complex biological neural network, and 
is now a computing system wildly applied in machine learning
\cite{haykin2009neural,krizhevsky2012imagenet,lecun2015deep}. 
Feedforward network with full connection between neighboring layers is one of 
the first models introduced \cite{mcculloch1943logical,rosenblatt1958perceptron}, 
and the algorithms evaluating and training them have been studied since then
\cite{rosenblatt1958perceptron,rosenblatt1962principles,minsky1969perceptrons,
werbos1974beyond}. Besides applications in image recognition 
\cite{krizhevsky2012imagenet}, natural language processing \cite{lecun2015deep},
cognitive science \cite{lake2015human}, and genomics \cite{alipanahi2015predicting},
neural network is also a powerful tool for function approximation
\cite{hornik1989multilayer,cybenko1989approximation,
jones1990constructive,carroll1989construction,liu2019neural}.
It is proved that functions in the form of multilayer feedforward network (MFN)
is dense in function spaces such as $C(I)$ and $L^2(I)$ ($I$ is unit
cube\footnote{The unit cube on $\mathbb{R}^n$ is defined as $[0, 1]^n$.}).
It is also easily shown that increasing the layers of MFN will enormously 
increase its function approximation capability. However, deep neural networks 
are difficult to train with gradient methods such as backpropagation due to 
gradient vanishing \cite{hochreiter1991untersuchungen,hochreiter2001gradient}. 
Here four-layer feedforward networks are chosen to avoid using special 
techniques such as ResNet \cite{he2016deep}.

Due to the ability of NN in function approximation, lots of efforts have been 
made to construct ODE/PDE solvers based on NN 
\cite{mall2014chebyshev,berg2017unified,raissi2017physics}. 
One of the major difficulties in such solvers is how to train a particular NN 
to satisfy the boundary condition accurately, since that the original form of 
NN as a trial function does not match boundary condition like trial functions 
in Galerkin methods. One strategy is the penalty method 
\cite{raissi2017physics,liu2019neural,wei2018machine}. The penalty method has 
been applied successfully in Burgers equation \cite{raissi2017physics}, Laplace 
equation \cite{liu2019neural}, and diffusion equations \cite{wei2018machine}, 
but only limited accuracy can be achieved. Another issue is how to evaluate the 
derivatives in equations, which need to be compatible with the NN-based solver. 
One option is the so-called automatic differentiation (AD) 
\cite{raissi2017physics}. AD evaluates the derivative with respect to input 
variable of any function defined by a computer program, and it is done by 
performing a non-standard interpretation of the program: the interpretation 
replaces the domain of the variables and redefines the semantic of the 
operators \cite{rall1981automatic}.

In our framework, we define a trial function which consists of a bulk term and 
a boundary term. The boundary term matches initial/boundary conditions, and the 
bulk term satisfies a reduced problem with relaxed boundary constrains,
respectively. The boundary term can be construct explicitly. We then define for 
the bulk term a loss function, which is actually the residual of the reduced 
problem. Such loss function does not involve any boundary conditions since the 
boundary conditions are relaxed in the reduced problem. Finally the bulk term, 
and therefore the trial function, is determined by minimizing the loss 
function. Machine learning technique is used for this minimization of the loss 
function. We refer to this new strategy as the constrained multilayer 
feedforward network (CMFN) method. With such novel strategy we will show that 
much higher accuracy can be achieve. It should also be pointed out that any 
method can be used to minimizing the loss function.

Before proceeding, we would like to clarify some terminology. In the language 
of the machine learning community, the trial function is usually called 
\emph{model}. The minimization of the loss function is actually a 
\emph{learning process}, during which the trial function \emph{learns} 
the correct data distribution of analytical solution. The minimization process 
is also a standard \emph{optimization} problem, and it is equivalent to 
\emph{training} in machine learning. Thus, the terminologies ``training'',
``optimization'', and ``minimization'' will be used interchangeably throughout 
this paper.

The paper is organized as follows. In section~2 we describe the framework in 
detail. Section~3 presents some numerical examples. Finally section~4 concludes 
the paper.

\section{Numerical Method}

To solve ODEs/PDEs numerically, one finds a function which satisfies the 
differential equations inside the domain and all initial/boundary conditions at 
(temporal/spatial) boundaries. That is, two parts of information need to be 
transferred into the numerical solver. For instance, in FVM the former part of 
information is transferred by flux reconstruction and the later part by 
operations on the boundary cells, respectively. In CMFN method, the former part 
of information is transferred by directly applying the differential operators 
with AD technique. The initial/boundary conditions are dealt with the boundary 
term in the trial function.

The CMFN method is based on the concept of MFN\@. MFN with $n$ layers can be 
defined as a computing algorithm as follows. The input layer as a vector is 
denoted by $y^{(1)}$, the output layer is denoted by $y^{(n)}$, and the hidden 
layers by $y^{(i)}$.  The output layer $y^{(n)}$ is computed by hidden layer 
$y^{(n-1)}$:
\begin{gather*}
    y^{(n)}_k = \sum_j \theta^{(n-1)}_{kj} y^{(n-1)}_j + \beta_k^{(n-1)},
\end{gather*}
and the hidden layers are computed recursively by:
\begin{gather*}
    \begin{cases}
	z^{(i+1)}_k = \sum_j \theta^{(i)}_{kj} y^{(i)}_j + \beta_k^{(i)} \\
	y^{(i+1)}_k = \phi(z^{(i+1)}_k)
    \end{cases}
    \quad i = 1, 2, \ldots, n-2.
\end{gather*}
Explicitly, a three-layer feedforward network
$N(x;\theta, \beta)$ is defined with superposition of activation function $\phi$ over 
linear transformation ($x=\{x_i\}_{n\times1}$ as input layer and 
$N(x; \theta, \beta) = \{N_k\}_{m\times1}$ as output layer):
\begin{gather}
    \label{eq:3layerMFN}
    N_k = \sum_j \theta_{kj}^{(2)}\phi(\sum_i \theta_{ji}^{(1)}x_i+\beta_j^{(1)})
    +\beta_k^{(2)},
\end{gather}
and a four-layer MFN is
\begin{gather}
    \label{eq:4layerMFN}
    N_l = 
        \sum_k \theta_{lk}^{(3)}\phi(
            \sum_j \theta_{kj}^{(2)}\phi(
	        \sum_i \theta_{ji}^{(1)}x_i+\beta_j^{(1)})
	       + \beta_k^{(2)}) 
	   + \beta_l^{(3)}.
\end{gather}
The parameters of the MFN are its weights $\theta = \{\theta_{ij}^{(k)}\}$ and 
biases $\beta=\{\beta_j^{(k)}\}$. It has been proved that the MFN 
\cref{eq:4layerMFN} with proper activation $\phi$ is dense in $C(I)$, namely 
the set of all continuous functions defined on unit cube 
\cite{cybenko1989approximation}. Therefore, for any continuous function $y(x)$ 
defined on a finite domain, one set of parameters $(\theta^*, \beta^*)$ can be 
found such that the corresponding network $N(x; \theta^*, \beta^*)$ is close 
enough to $y(x)$, i.e.\ the norm $\Vert y(x) - N(x; \theta^*, \beta^*)\Vert $ 
could be sufficiently small. Similar conclusion holds for $y \in L^2(I)$, which 
is $\int_I |y(x)|^2 \dd{x} <\infty$. Such properties guarantee that an optimal 
set of $N(x; \theta^*, \beta^*)$ exists with the corresponding MFN being a good 
numerical approximation of the solution.

A well-posed ODE/PDE with Dirichlet boundary condition can be written as 
\begin{equation}
    \begin{cases}
	\mathcal{L}u = f &\qqtext{in} \Omega, \\
	\mathcal{B}u = g &\qqtext{in} \partial \Omega. 
    \end{cases}
    \label{eq:general-DE}
\end{equation}
In CMFN we define a model function:
\begin{equation}
    \label{eq:CMFN-basic}
    \hat{y}(x; \theta, \beta) = G(x) + \tilde{N}(x;\theta,\beta) \equiv G(x) + w(x)\cdot N(x;\theta,\beta).
\end{equation}
As the boundary operator $\mathcal{B}$ is linear and algebraic, we choose
the two terms $G(x)$ and $w(x)$ in \cref{eq:CMFN-basic} such that
\begin{enumerate}
    \item $\mathcal{B}G=g$ when $x \in \partial \Omega $,
    \item $\mathcal{B}\tilde{N}\to0$ when $x\to\partial\Omega$.
\end{enumerate}
$G(x)$ is a boundary term which is a pre-defined function, $\tilde{N}(x)$ is 
bulk term, and $N(x)$ is the unknown part which is approximated by a neural 
network. The original problem with respect to $u(x)$ in \cref{eq:CMFN-basic} is 
reduced to solving a new differential equation with respect to $\tilde{N}(x)$. 
The the new equation is defined as the \emph{reduced equation}, and the unknown 
part $N(x)$ separated from bulk term is called \emph{reduced solution}.

In \cref{eq:CMFN-basic}, as long as a pre-defined weight $w(x)$ is a bounded 
continuous function, the bulk term $\tilde{N}(x)$ is continuous and bounded 
according to \cref{eq:4layerMFN}. The bulk term is further written as 
$\tilde{N}(x) = w(x)\cdot N(x)$, where the pre-defined weight $w(x)$ satisfies:
\begin{enumerate}
    \item $w(x)\to0$ as $x\to\partial\Omega$ (vanishing on domain boundary),
    \item for all $x^* \in\Omega$, $w(x^*)\neq 0$ (non-vanishing within domain).
\end{enumerate}
Now the boundary conditions are automatically satisfied, or we say, the 
boundary conditions are \emph{relaxed} in the reduced equation. Once the 
reduced equation is determined, the loss function can be constructed without 
considering the original boundary conditions. 

By substituting the trial function $\hat{y}(x) = G(x) + w(x)\cdot N(x)$ into 
\cref{eq:general-DE}, the original problem is converted to its reduced 
equation:
\begin{equation}
    \label{eq:reduced-general}
    \mathcal{L}[G + w\cdot N] = 
        \tilde{\mathcal{L}}[N] = f \quad N \in C(\mathbb{R}^n).
\end{equation}
One may think that there could be \emph{multiple} solutions to reduced equation 
since it has no boundary conditions. However, while the original problem 
\cref{eq:general-DE} has unique solution $y^*$, the reduced one 
\cref{eq:reduced-general} should also have a \emph{unique} solution 
$(y^*-G)/w$. The paradox indicates that, among all solutions to 
\cref{eq:reduced-general}, there exists a unique solution satisfies: 
$\mathcal{B} \tilde{N} \to 0$ while $x\to\partial\Omega$. We do not have a 
rigorous proof for this statement yet, but as supported by the examples shown 
later, a unique solution can always be obtained.

After construction of trial function $\hat{y}$, the loss function towards which
optimization is done is defined by residual
$\mathcal{R}N = \tilde{\mathcal{L}}N-f$ of \cref{eq:general-DE}:
\begin{equation}
    \label{eq:loss-general}
    L = 
    \int_\Omega \left\langle(\mathcal{R}N)(x),(\mathcal{R}N)(x)\right\rangle  
    \dd{x} = 
    \sum_{x^*\in T(\Omega)}  \Vert(\mathcal{R}N)(x^*)\Vert^2 ,
\end{equation}
where $T(\Omega)$ is the training set containing points selected from domain 
$\Omega$. The operator $\mathcal{R}$ is defined by AD instead of manually 
working it out. This not only saves the researches from laborious job 
\cite{berg2017unified}, but also produces robust and reliable code 
\cite{baydin2015automatic,rall1981automatic}. There are successful AD 
implements on nearly all programming platforms \cite{baydin2015automatic}.  In 
this work, reverse mode AD application programming interface on TensorFlow 
\cite{tensorflow2015-whitepaper} is called. Since AD solves the problem of much 
too complicated differential problem, high order differential operator 
$\mathcal{L}$ are solved in this work without extra efforts.

The final stage of the framework solving ODE/PDE is optimization during which 
loss function $L$ defined by \cref{eq:loss-general} is minimized with respect 
to its free parameters. In cases where MFN $N(x;\theta, \beta)$ is the reduced 
solution, its weights and biases $(\theta, \beta)$ are trained to minimize 
$L=L(\theta, \beta)$. In this work, the optimization is done by second order 
method L-BFGS \cite{Liu1989on}, instead of SGD, the most popular choice 
in building machine learning models \cite{lecun2015deep}. The second order 
method is not always robust in general machine learning problems, but it serves 
well in ODE/PDE solver according to our observation. All numerical examples 
presented in later section in this paper is trained by second order method 
L-BFGS which greatly improves training efficiency. The training process 
requires large amount of computational resource which used to be an obstacle in 
development of machine learning \cite{minsky1969perceptrons}. Parallelism and 
heterogeneous computing throw lights on the problem, and the model in this work 
is defined and trained on TensorFlow \cite{tensorflow2015-whitepaper}.

\section{Examples and Discussion}

\begin{figure}[hptb]
    \centering
    \subfloat[][Model $\hat y(x)$\label{fig:1d-integral-y}]{%
	\includegraphics[width=0.45\textwidth]{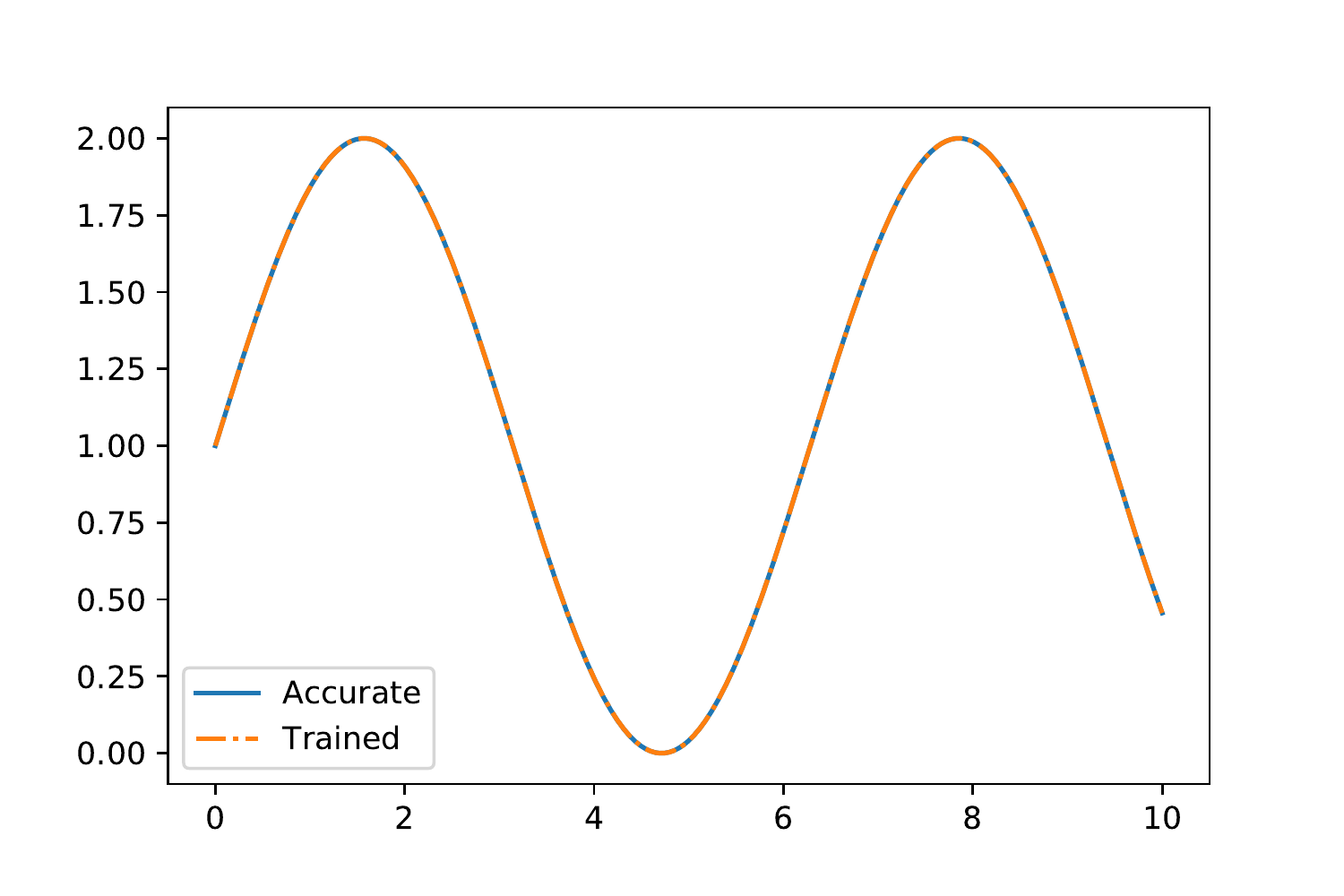}
    }
    \subfloat[][Error distribution of Model\label{fig:1d-integral-y-dev}]{%
	\includegraphics[width=0.45\textwidth]{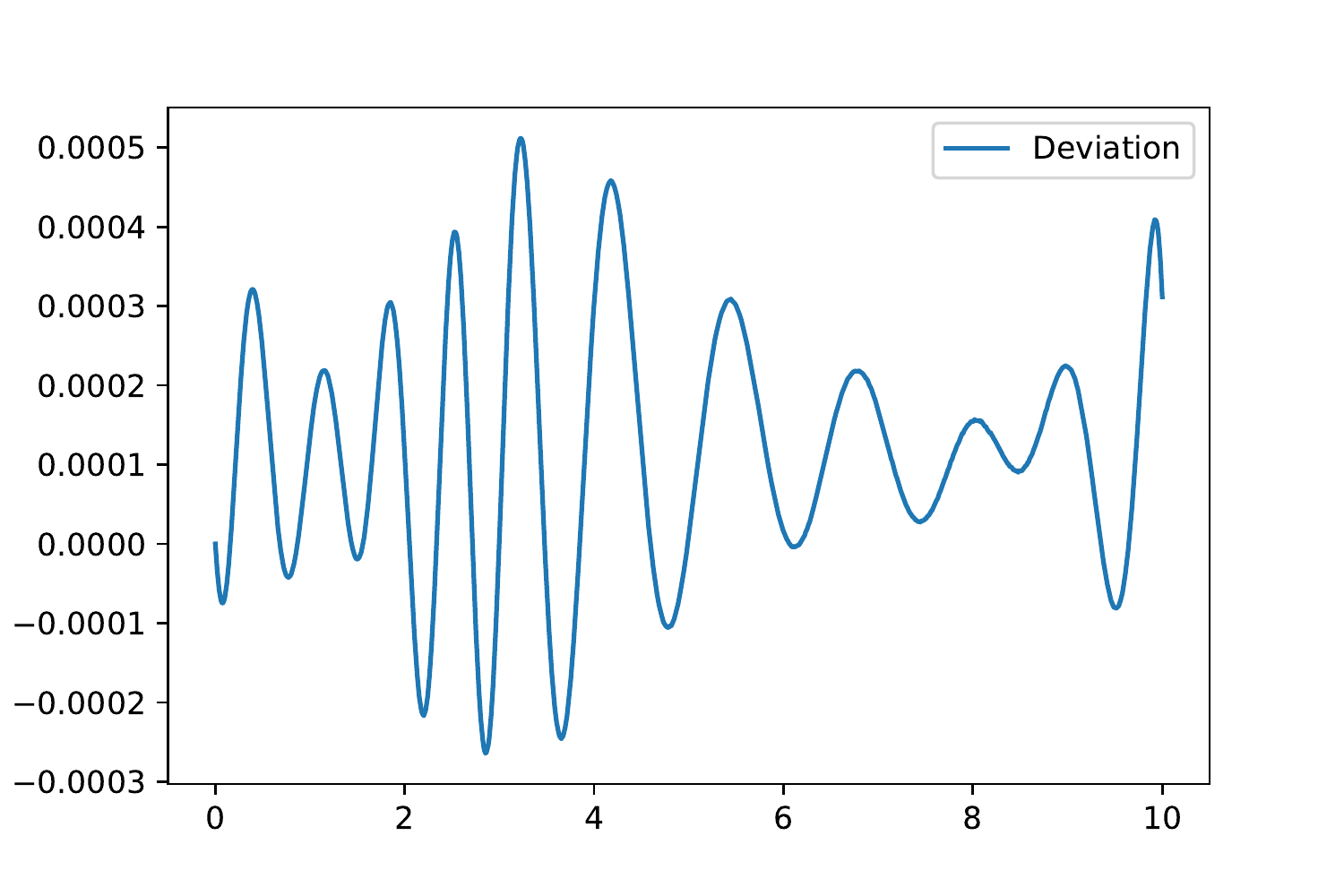}
    }\\
    \subfloat[][Bulk term $w(x)\cdot N(x)$\label{fig:1d-integral-core}]{%
	\includegraphics[width=0.45\textwidth]{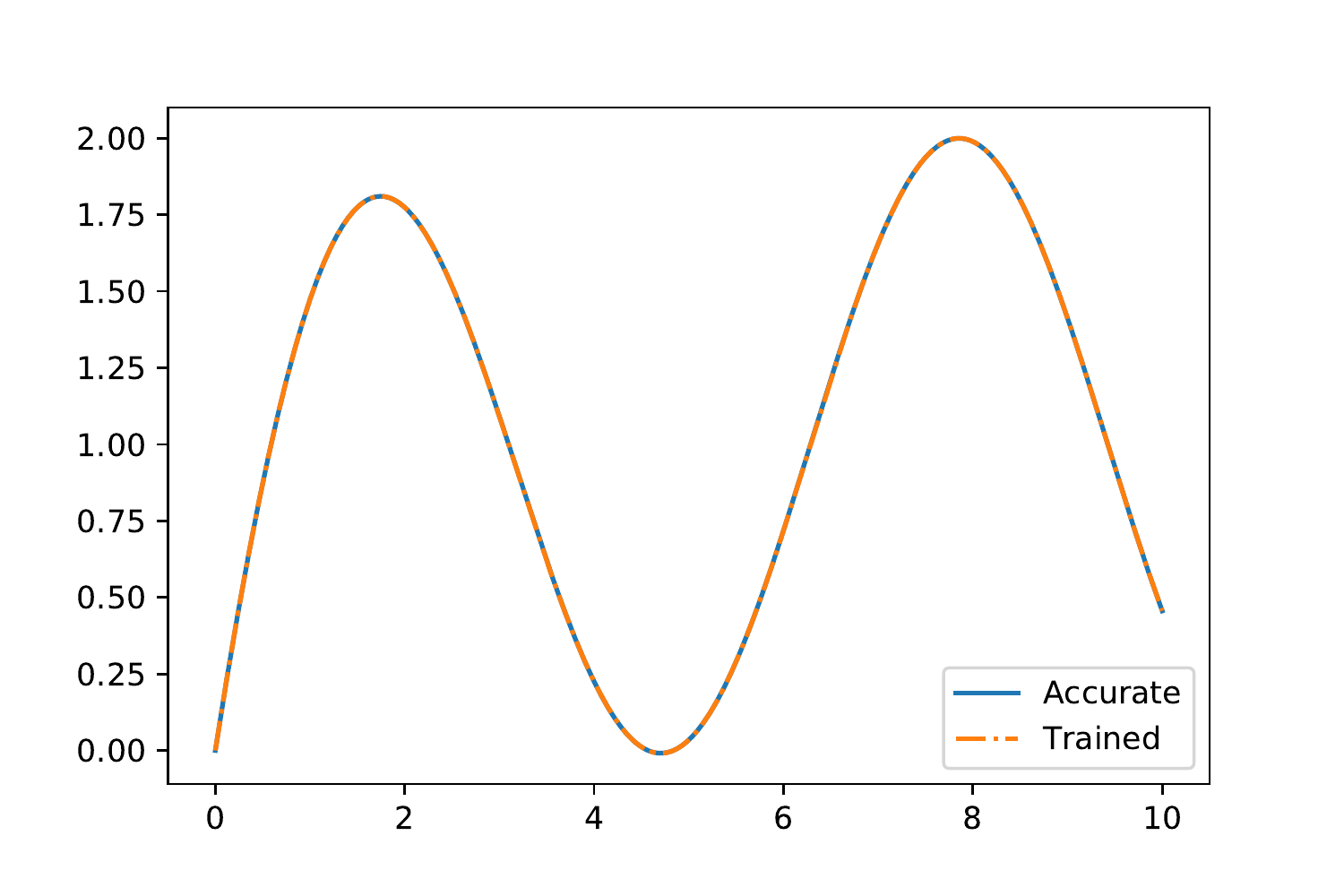}
    }
    \subfloat[][Error distribution of bulk
	term\label{fig:1d-integral-core-dev}]{%
	\includegraphics[width=0.45\textwidth]{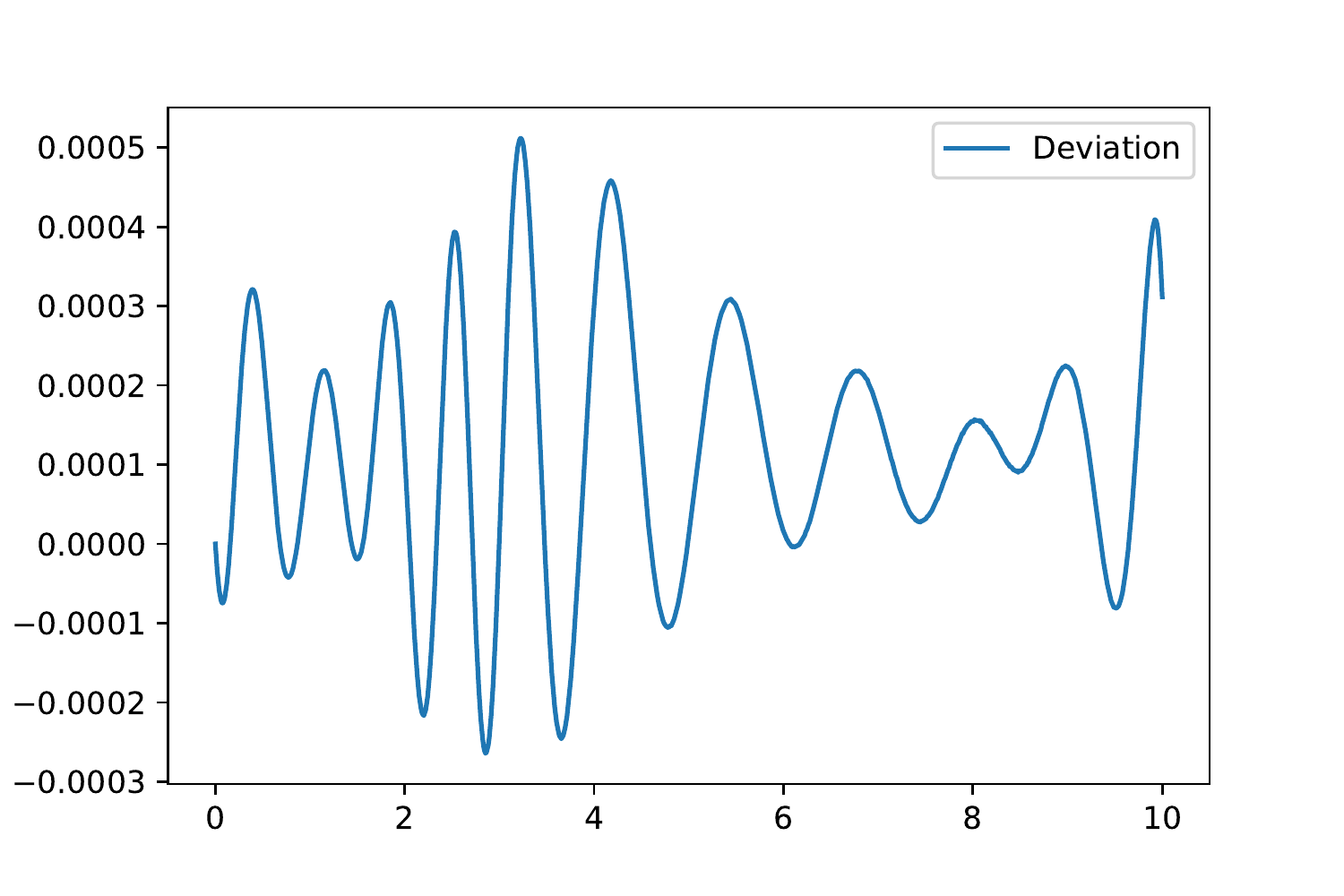}
    }\\
    \subfloat[][Reduced solution $N(x)$\label{fig:1d-integral-N}]{%
	\includegraphics[width=0.45\textwidth]{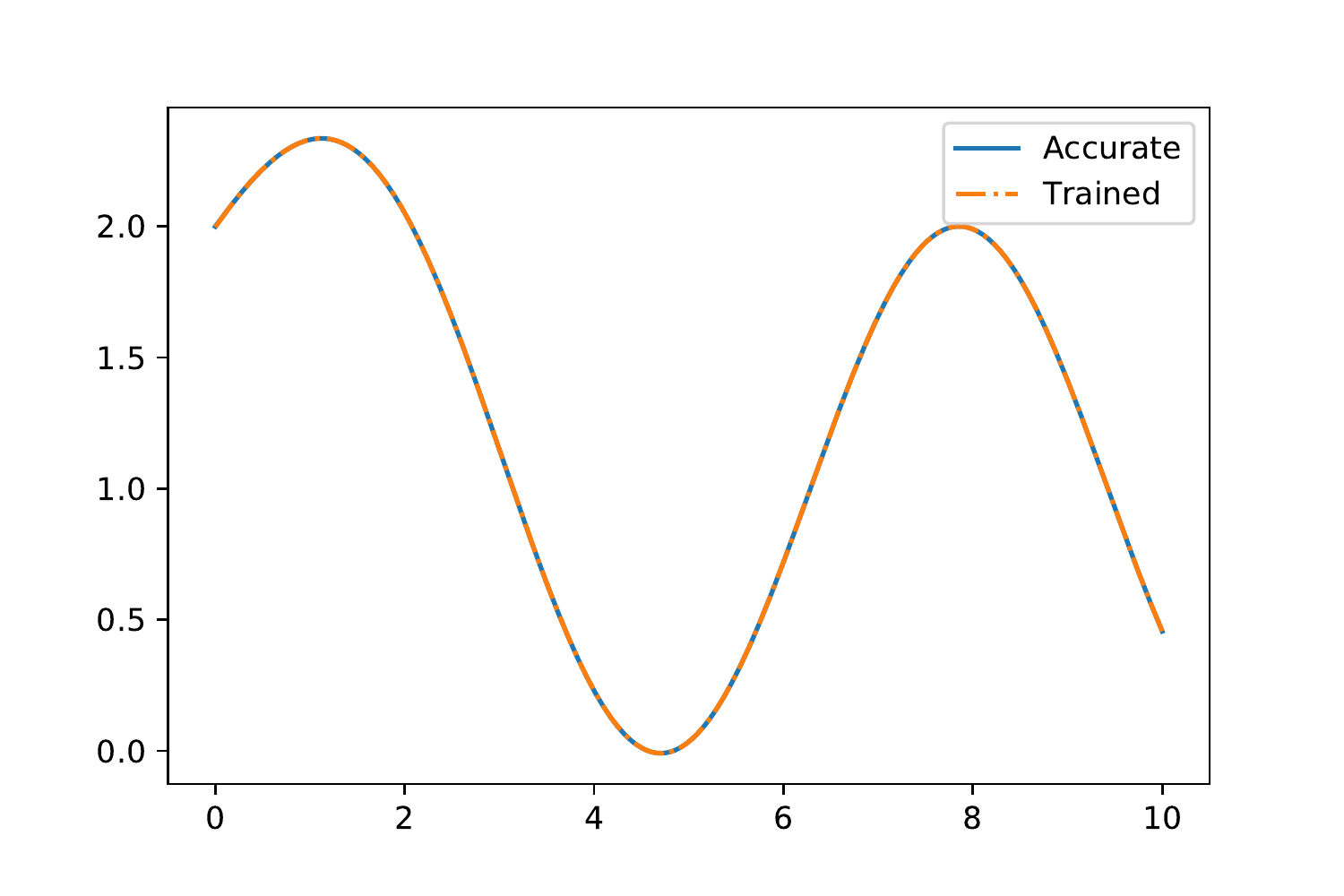}
    }
    \subfloat[][Error distribution of reduced solution
	\label{fig:1d-integral-N-dev}]{%
	\includegraphics[width=0.45\textwidth]{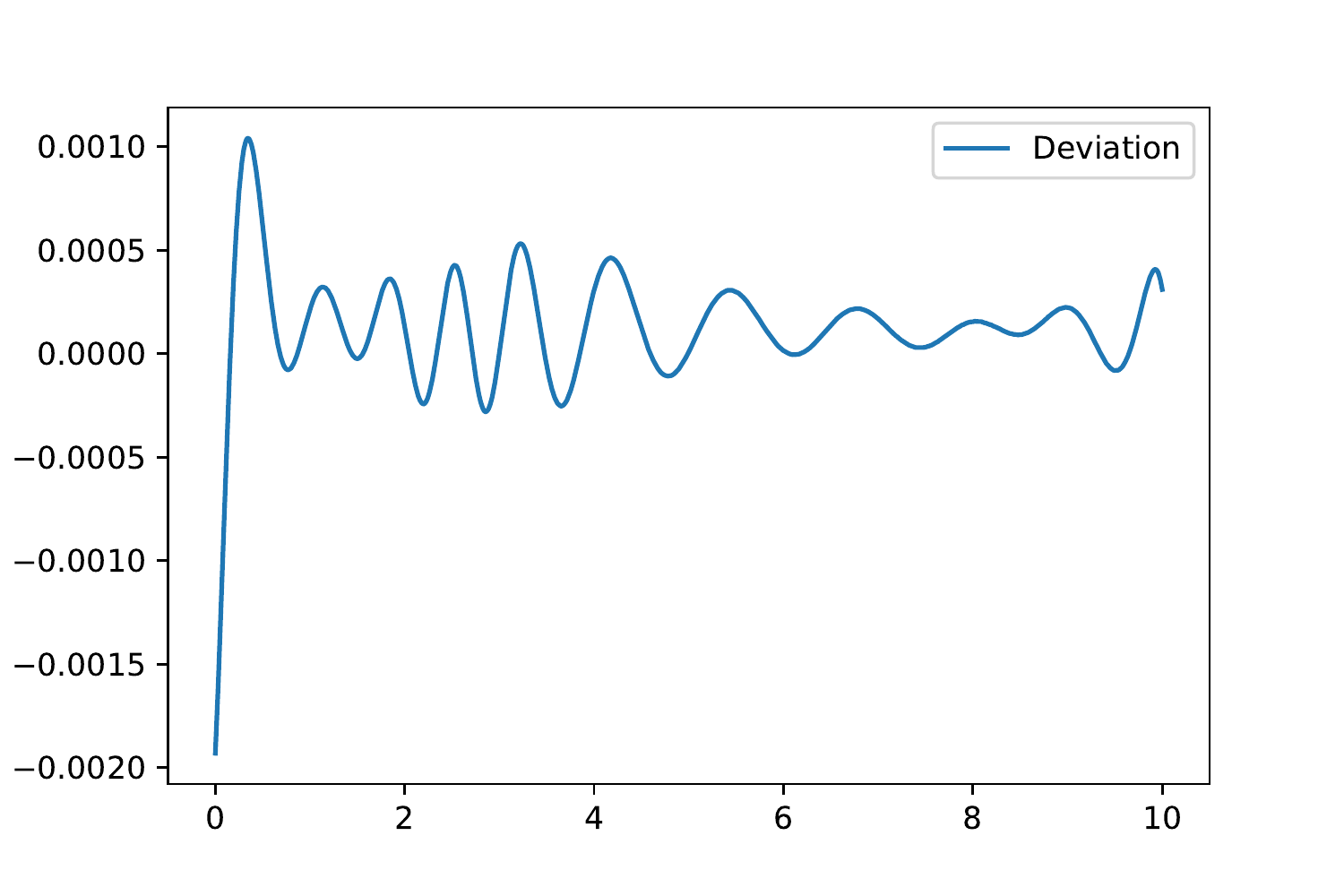}
    }
    \caption{Results of definite integral problem in \cref{eq:1d-integral} with 
	$\mathrm{Error} = 1.96\times10^{-4}$}
    \label{fig:1d-integral}
\end{figure}

The first example on ODE solving is a definite integral problem as 
illustration:
\begin{equation}
    \label{eq:1d-integral}
    \begin{cases} 
	y'(x) = \cos x \\
	y(x=0) = y_0 = 1 
    \end{cases}.
\end{equation}
The analytical solution is simply integration of R.H.S. of \cref{eq:1d-integral}:
$y(x)=1 + \sin x$. In order to find a numerical solution in domain $[0, 10]$. The trial
function is defined as
\begin{equation}
    \label{eq:1d-integral-model}
    \hat{y}(x;\theta,\beta) = y_0 e^{-x} + (1-e^{-x})N(x;\theta,\beta).
\end{equation}
It is easily verified that requirements for $G(x)$ and $\tilde{N}(x)$ are all 
satisfied. The network $N(x;\theta, \beta)$ is set as a four-layer network with 
$20$ neurons in each hidden layer. Loss function is defined as
\begin{equation}
    \label{eq:1d-interal-loss}
    L = \int_0^{10} \Vert \hat{y}' - \cos x\Vert^2 \dd{x} =
    \sum_{i=1}^{1000} | \hat{y}'(x_i) - \cos x_i |^2,
\end{equation}
with $\{x_i\}_{i=1,2,\ldots,1000}$ are uniformly selected in interval $[0, 10]$.
The loss function is minimized by L-BFGS method, the result is illustrated in
\cref{fig:1d-integral}.

The boundary condition is well maintained in this example. 
\Cref{fig:1d-integral-y,fig:1d-integral-core,fig:1d-integral-N} shows that the 
trained MFN provides a correct reduced solution, as well as a good match in 
bulk term and the whole model. However, error distributions shown in
\cref{fig:1d-integral-y-dev,fig:1d-integral-core-dev,fig:1d-integral-N-dev}
reveal that the reduced solution actually has the lowest accuracy, especially 
on domain boundary, the pre-defined constrains in $G(x)$ and $w(x)$ reduce 
error in model, as well as turn the training process into a standard 
non-constrained optimization problem.

Another important property shown in \cref{fig:1d-integral-y-dev} is that
the solution learned by CMFN deviates from analytical solution randomly, where
traditional numerical methods usually have increasing errors along iterations.
This property is natural for CMFN method since all data points are equal to the
learning process, while in iteration methods error grows accumulatively. A more 
accurate measure of error in solution provided by \cref{eq:1d-integral-model} 
is calculating the $L^2$-norm of error distribution:
\begin{equation}
    \mathrm{Error} = 
        \sqrt{\frac{1}{10} \int_0^{10} |y(x) - \hat{y}(x)|^2 \dd{x}.}
\end{equation}
The previous example has an average error of $1.96\times10^{-4}$. The authors 
also tried other choices on defining boundary term $G$ and weight function $w$: 
$G^*(x) \equiv 1$ and  $w^*(x) = x$. Numerical test shows that definitions in
\cref{eq:1d-integral-model} is more optimized than their alternatives. $G^*(x)$
instead of original $G(x)$ in \cref{eq:1d-integral-model} roughly doubles the 
average error, and $w^*(x)$ instead the original weight even reduces the 
accuracy for one order of magnitude.

As is mentioned above, the reduced equation with relaxed boundary condition 
should have unique solution instead of multiple solution under the premise that 
the solution is bounded. The reduced equation for \cref{eq:1d-integral} is
\begin{equation}
    \label{eq:1d-reduced}
    (1-e^{-x}) N' + e^{-x} N - y_0e^{-x} - \cos x = 0,
\end{equation}
and the general solution to \cref{eq:1d-reduced} is
\begin{equation}
    \label{eq:1d-reduced-solution}
    N(x) = \frac{C - y_0e^{-x}}{1-e^{-x}} + \frac{\sin x}{1-e^{-x}} \quad 
        C\in\mathbb{R}.
\end{equation}
At $x=1$, the first part in R.H.S. of \cref{eq:1d-reduced-solution} is 
unbounded if $C\neq y_0$; if $N(x)$ is assumed to be a bounded function on
$[0, 1]$, then in \cref{eq:1d-reduced-solution} there exists unique solution 
which is a proper reduced solution to problem in \cref{eq:1d-integral}.

The CMFN treats the problems defined by \cref{eq:general-DE} equally, so any 
initial condition problem is solved with similar process and accuracy. The 
following presents solution of boundary value problem (BVP) of a second order 
ODE\@. It is boundary layer problem \cite{prandtl1938zur} reduced to 1D\@:
\begin{equation}
    \label{eq:bl}
    u u' = \nu u'' \quad \quad u(0) = 1 \quad u(1) = 0.
\end{equation}
The problem has analytical solution:
\begin{equation}
    \label{eq:bl-solution}
    u(x) = \frac{2C}{1+\exp\left(\frac{x-1}{\nu}\cdot C\right)} - C\quad C > 1,
\end{equation}
with $C$ as a constant determined by algebraic equation:
\begin{equation*}
    1-\frac{2}{1+C} = \exp\left(-\frac{C}{\nu}\right).
\end{equation*}
We have $C \approx 1.2$ when $\nu=0.5$, and $C$ tends to $1$ rapidly as $\nu$ 
decreases to $\nu=0$.

The boundary term of model is constructed according to boundary conditions in 
\cref{eq:bl} where linear function $f(x)=1-x$ is sufficient in matching values 
of the two boundary points. The weight function is constructed by polynomial 
such that $w(0)=0$ and $w(1)=0$. Finally, the trial function is defined as
\begin{equation}
    \label{eq:bl-model}
    \hat{u}(x;\theta,\beta) = (1-x) + x(1-x)\cdot N(x;\theta, \beta),
\end{equation}
and loss function is defined similar to \cref{eq:1d-integral-model}. The model 
is trained with $100$ data points uniformly distributed in $[0, 1]$.

\begin{figure}[hpbt]
    \centering
    \subfloat[][Model $\hat{y}(x)$
	\label{fig:bl-y-dev}]{%
	\includegraphics[width=0.45\textwidth]{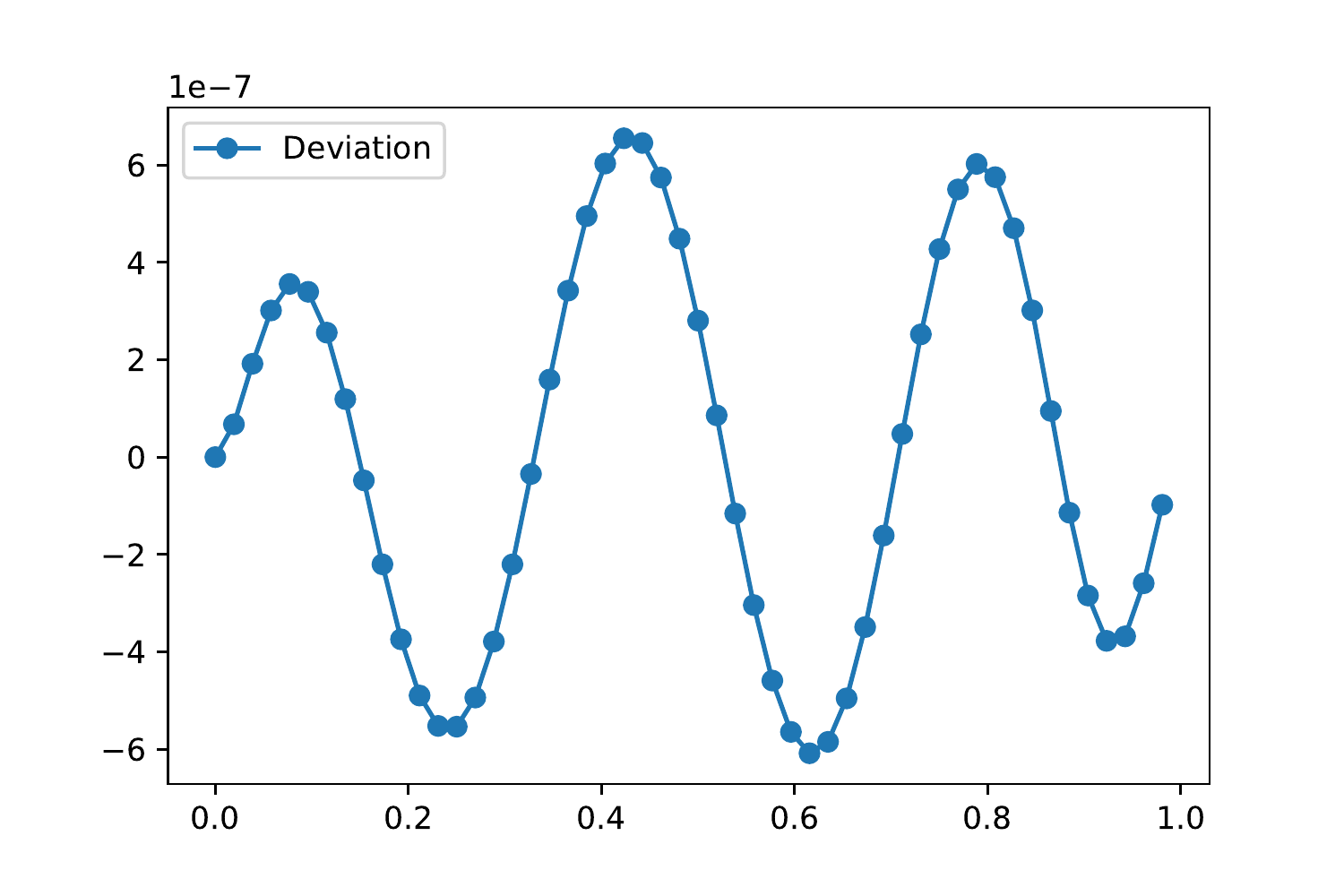}
    }
    \subfloat[][Reduced solution $N(x)$
	\label{fig:bl-N-dev}]{%
	\includegraphics[width=0.45\textwidth]{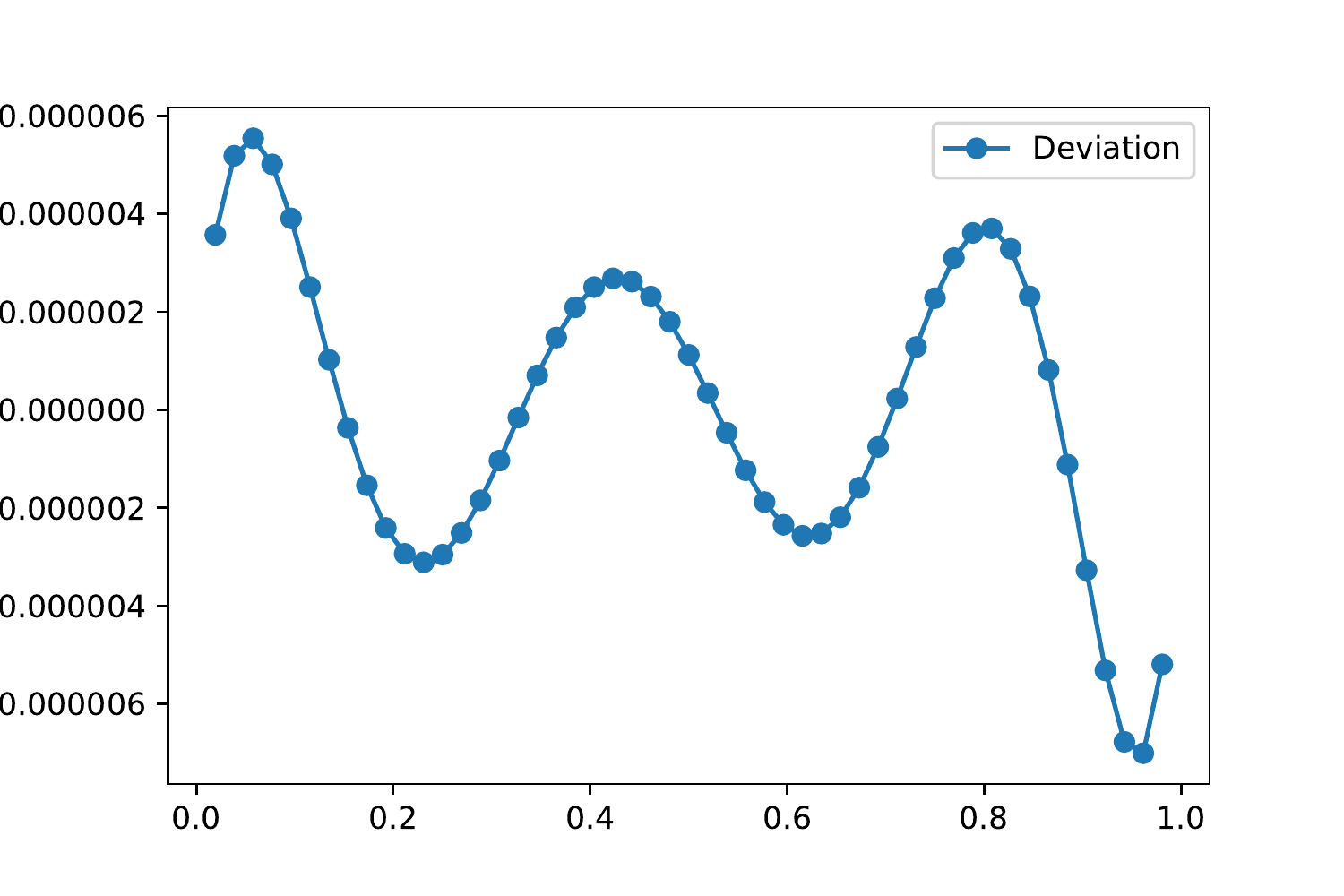}
    }
    \caption{Error distributions of BVP with average error $< 10^{-5}$}
    \label{fig:bl-dev}
\end{figure}

\begin{figure}[hptb]
    \centering
    \subfloat[][1st order derivative\label{fig:bl-diff1}]{%
	\includegraphics[width=0.3\textwidth]{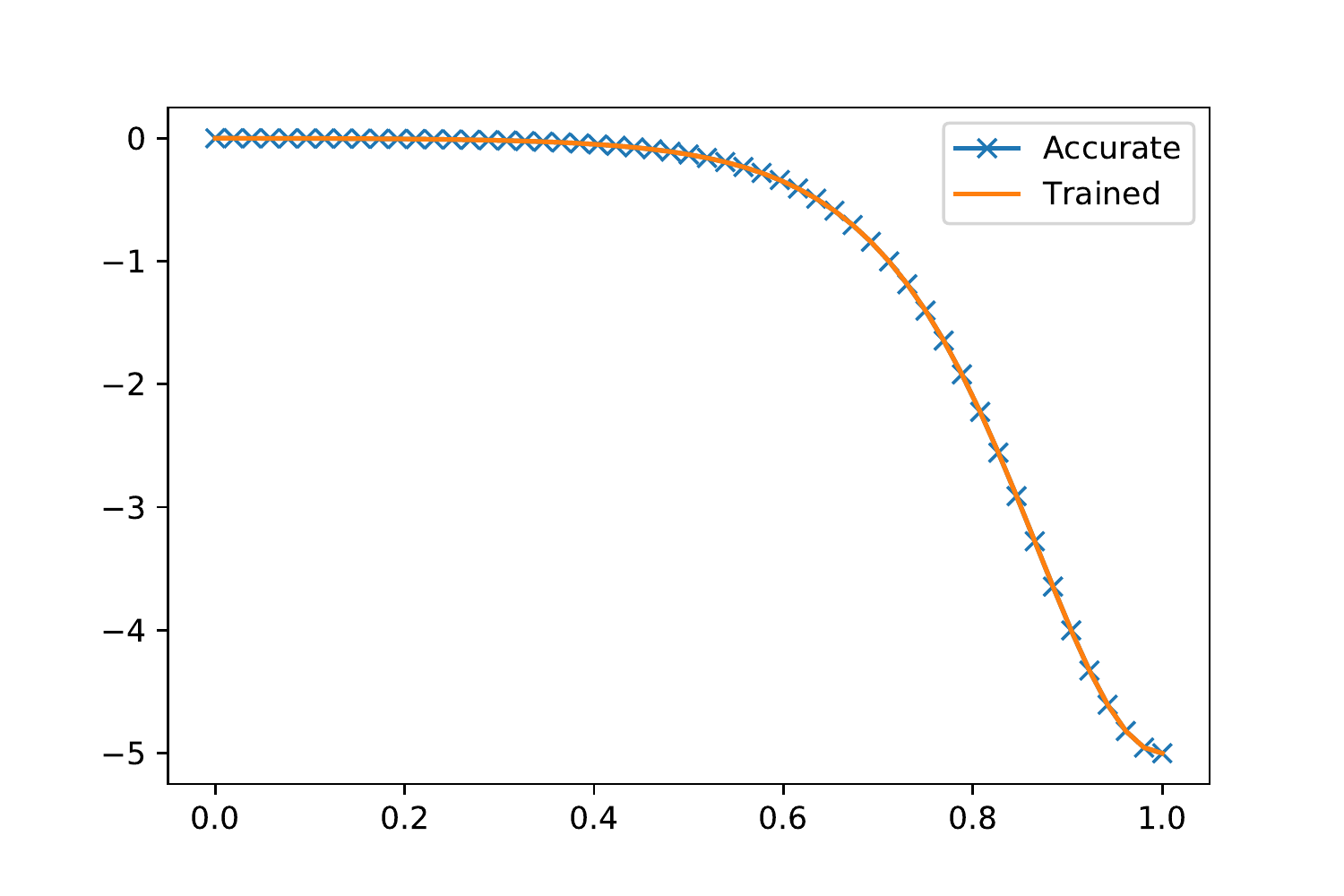}
    }
    \subfloat[][2nd order derivative]{%
	\includegraphics[width=0.3\textwidth]{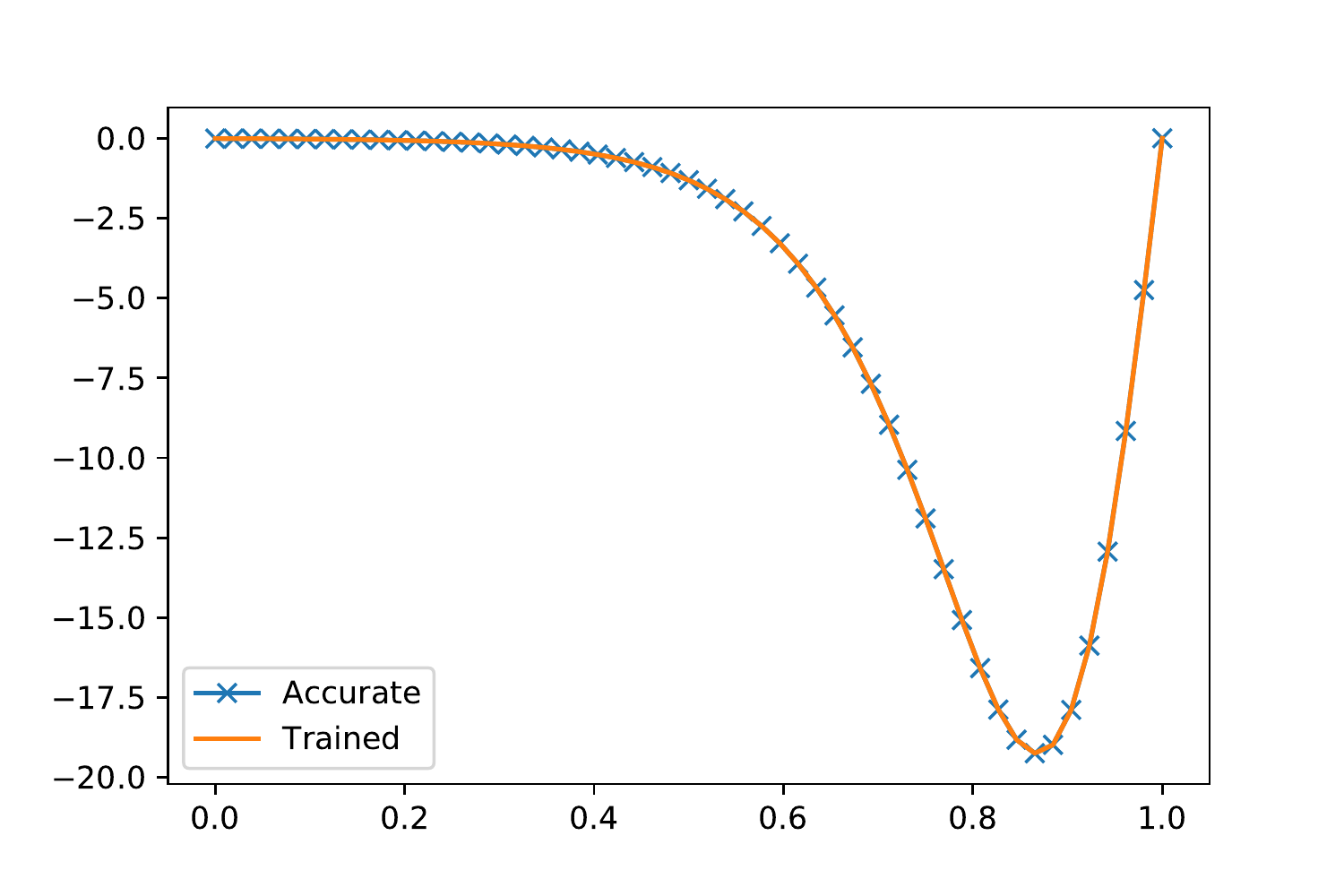}
	\label{fig:bl-diff2}
    }
    \subfloat[][3rd order derivative]{%
	\includegraphics[width=0.3\textwidth]{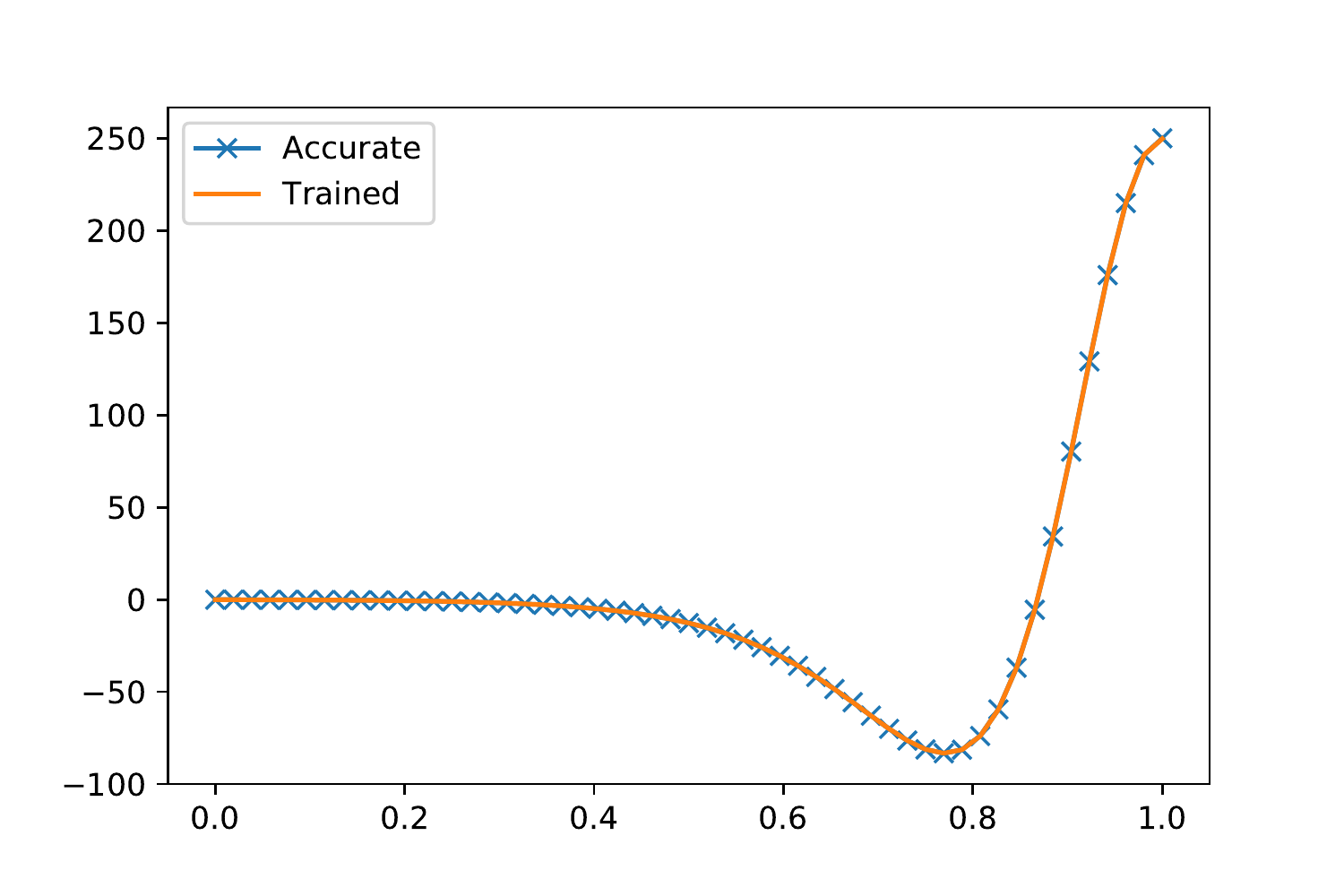}
	\label{fig:bl-diff3}
    }\\
    \subfloat[][4th order derivative]{%
	\includegraphics[width=0.3\textwidth]{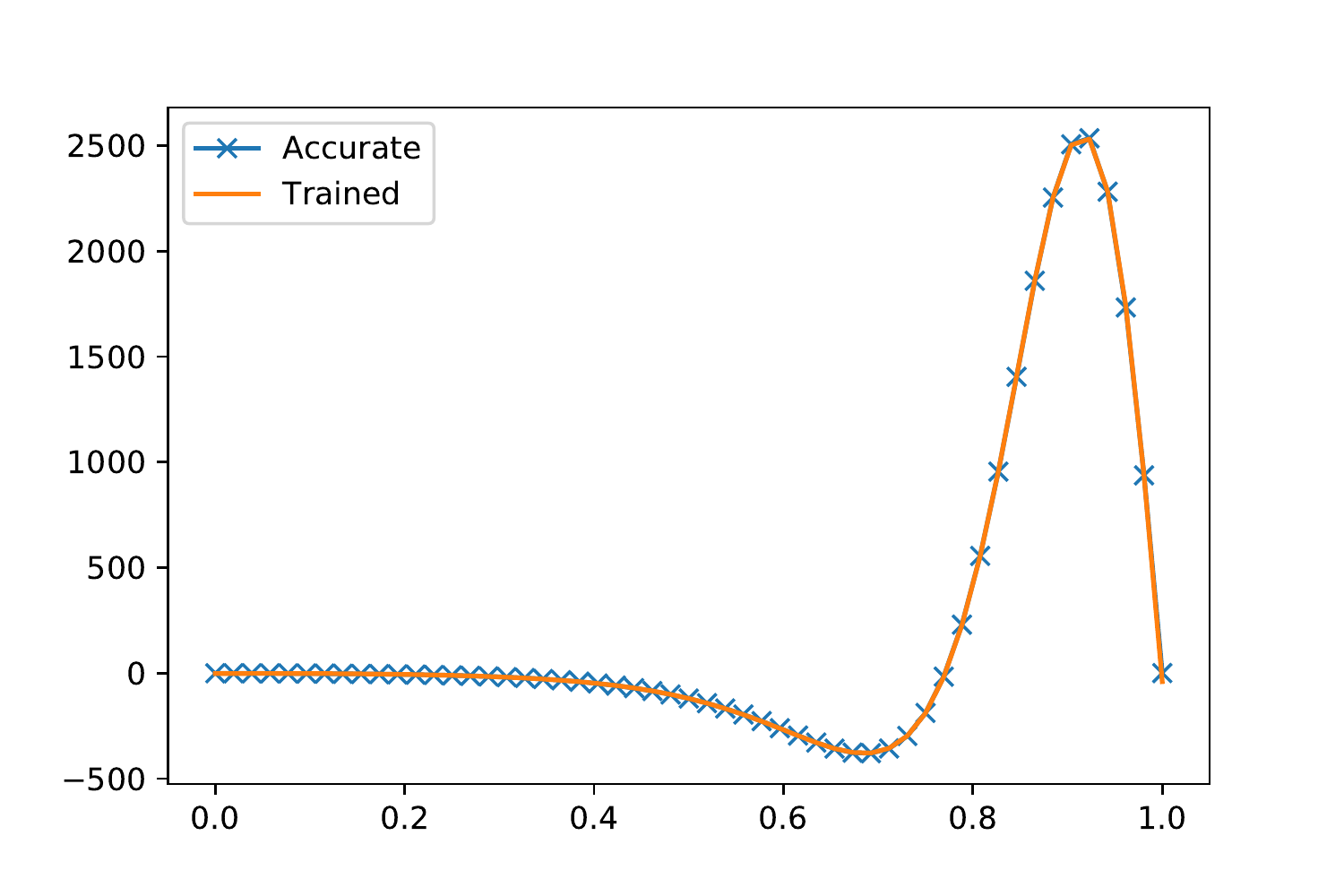}
	\label{fig:bl-diff4}
    }
    \subfloat[][5th order derivative]{%
	\includegraphics[width=0.3\textwidth]{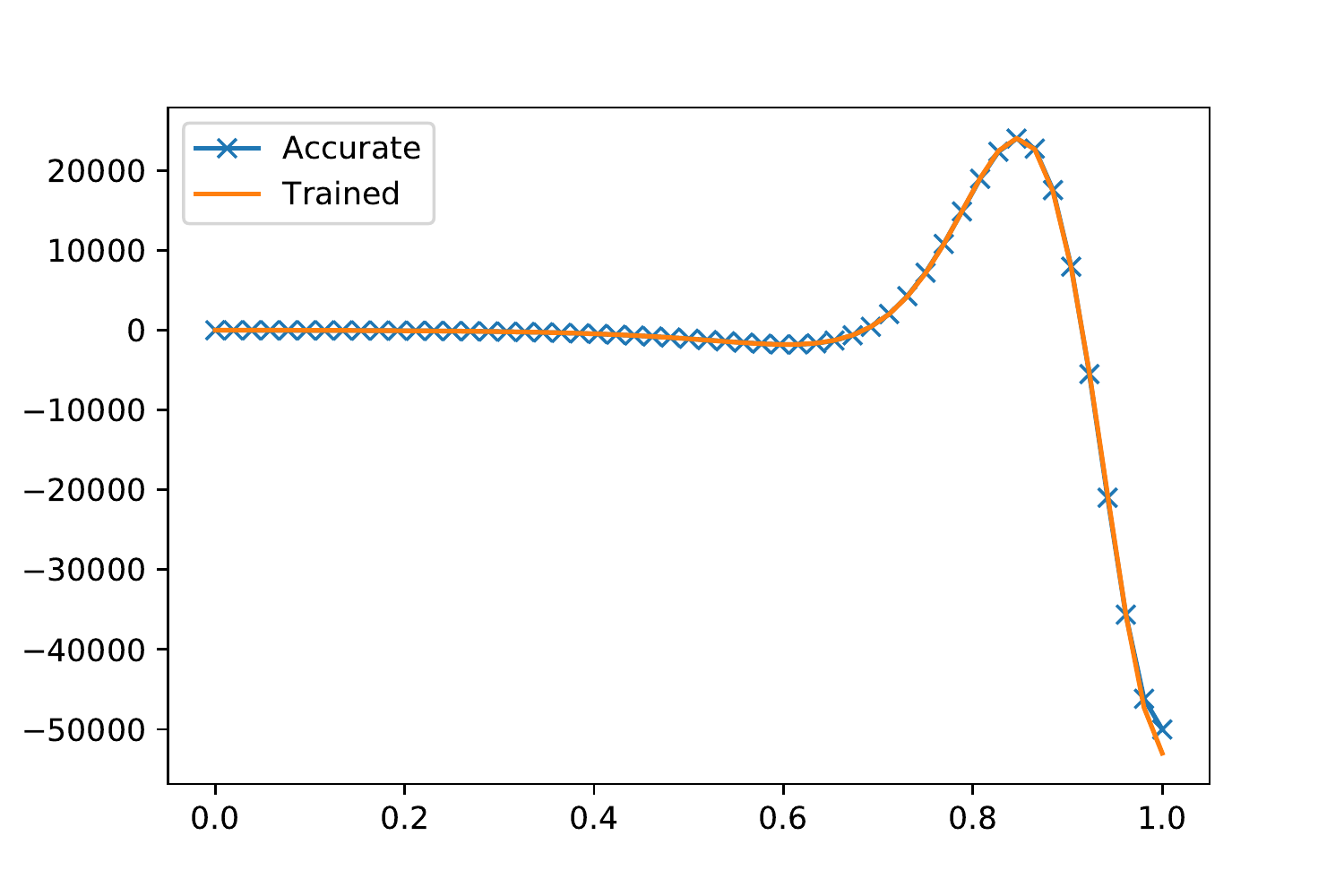}
	\label{fig:bl-diff5}
    }
    \subfloat[][6th order derivative]{%
	\includegraphics[width=0.3\textwidth]{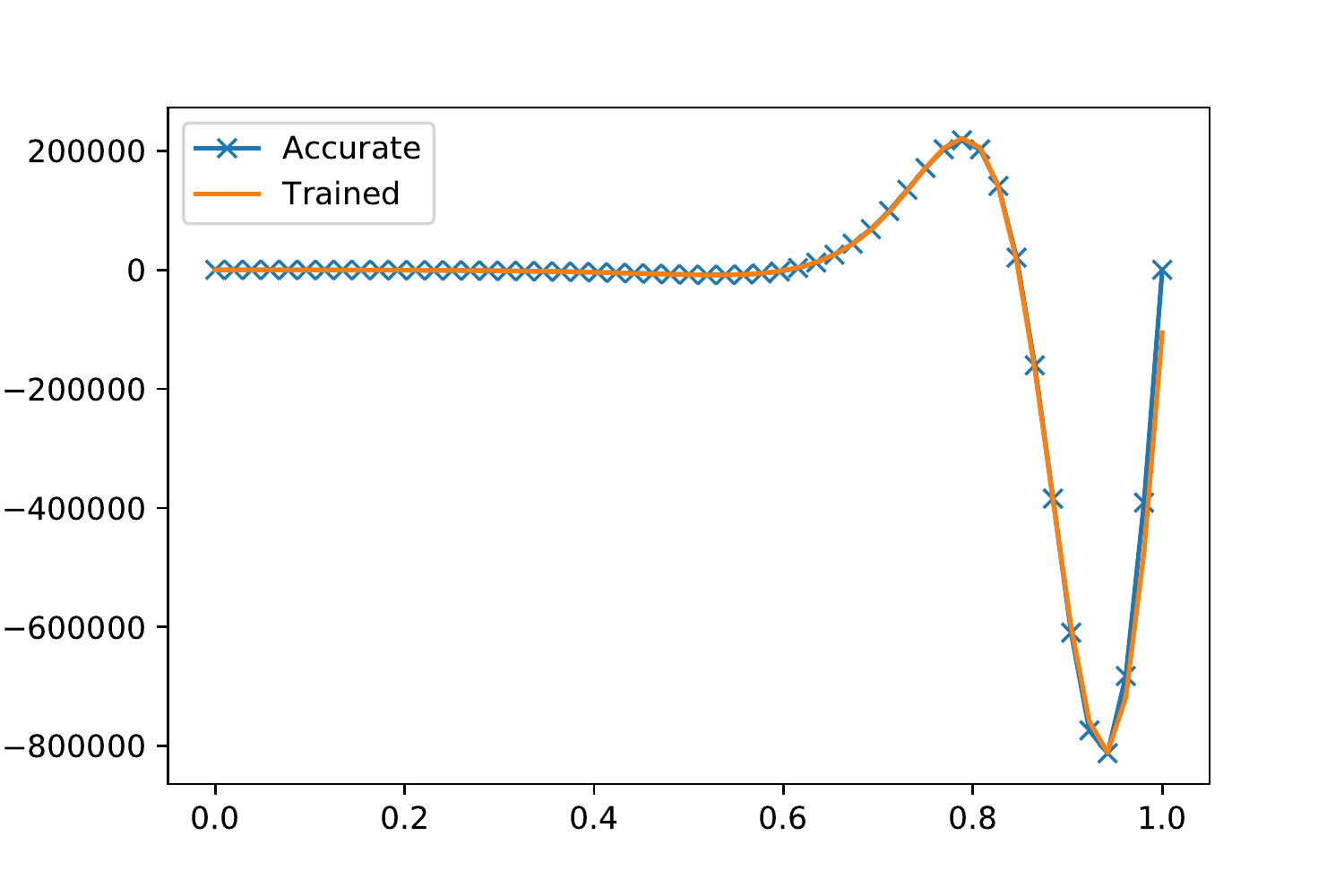}
	\label{fig:bl-diff6}
    }\\
    \subfloat[][7th order derivative]{%
	\includegraphics[width=0.3\textwidth]{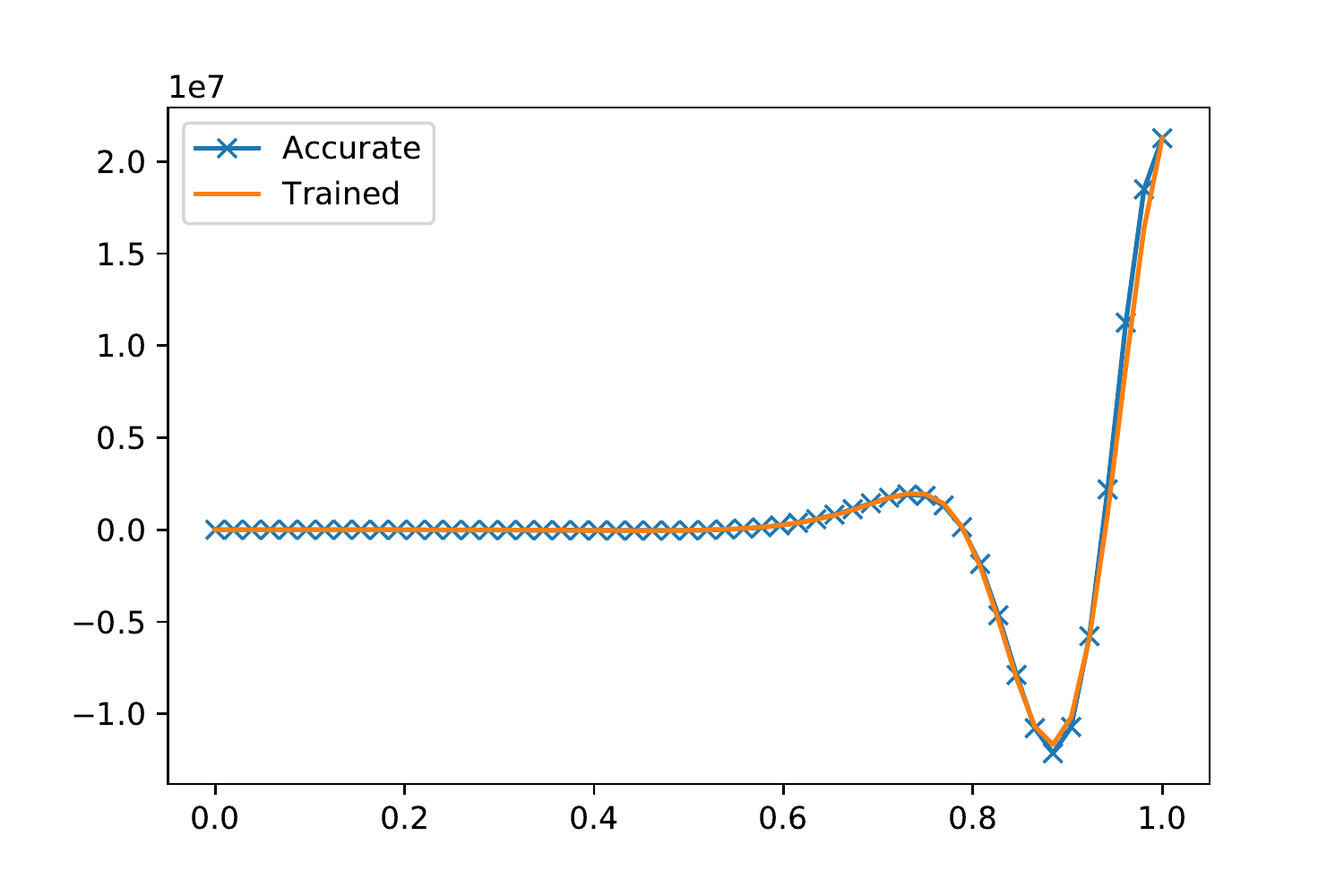}
	\label{fig:bl-diff7}
    }
    \subfloat[][8th order derivative]{%
	\includegraphics[width=0.3\textwidth]{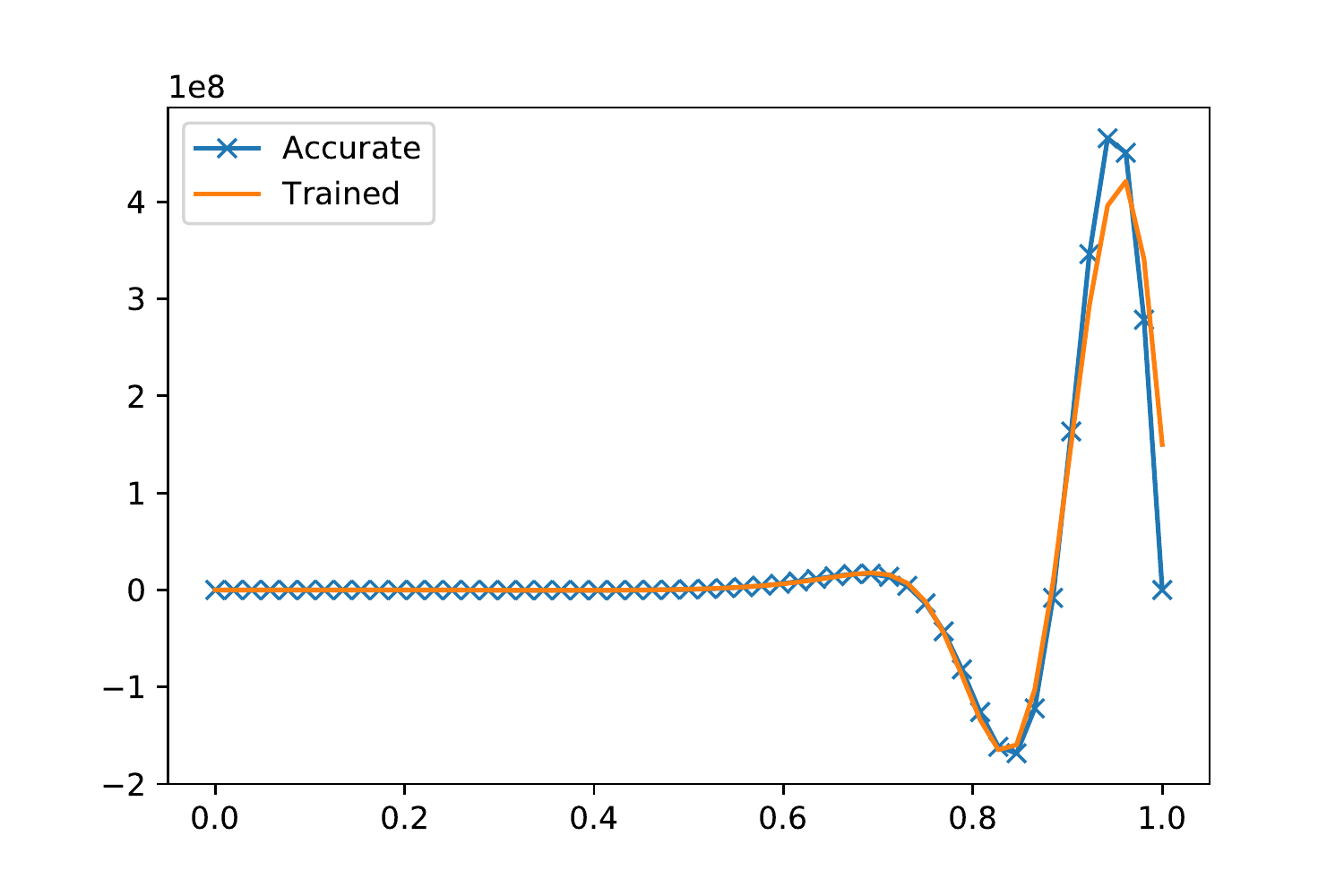}
	\label{fig:bl-diff8}
    }
    \caption{Differentiation of learned solution of \cref{eq:bl}}
\end{figure}

\Cref{fig:bl-dev} shows that BVP is solved with high numerical accuracy, and in 
\crefrange{fig:bl-diff1}{fig:bl-diff8} differential property of the learned 
solution is studied. CMFN solution to the problem is not only an accurate 
numerical approximation to the analytical solution, but its first to eighth 
order derivatives are also all accurate numerical approximations to their 
analytical counterpart. This is very hard to achieve with commonly used 
numerical methods enumerated in section~1.

Another important observation is that once weight is defined by polynomials,
there is a sole proper choice based on initial/boundary condition: for 1D 
Dirichlet boundary contion at $x=a$, there would be a factor $(x-a)$ in the 
weight function. For example, weight in \cref{eq:bl-model} should be defined as 
$w(x) = x(1-x)$, which is Hermite interpolation based on boundary condition: 
$w(1) = w(0) = 0$; if weight is defined as $w^*(x)=x^2(1-x)$ instead, there 
will be a huge reduction in accuracy. Since the term $x^2$ in $w^*(x)$ not only 
vanishes itself at $x=0$, but also has zero first order derivative there; as a 
result, the boundary term $G(x)$ unexpectedly dominates both function value and 
first order derivative at $x=0$.

Two dimensional problem is tested with heat conduction problem 
(Laplace equation):
\begin{equation}
    \label{eq:laplace}
    \nabla^2 T(x, y) = 0 \quad\quad x\in[0, 1] \quad y\in[0, 1]
\end{equation}
The boundary condition of Dirichlet problem is:
\begin{equation}
    \label{eq:dirichlet-bd}
    T(0, y) = T(1, y) = T(x, 0) = 0 \quad T(x, 1) = \sin \pi x.
\end{equation}
The problem has analytical solution:
\begin{equation}
    \label{eq:heat-solution}
    T(x, y) = \frac{\sin\pi x\sinh\pi y}{\sinh\pi},
\end{equation}
so the error of numerical solution $\hat{T}(x, y)$ is evaluated as:
\begin{equation}
    \label{eq:heat-error}
    \mathrm{Error} = \sqrt{\frac{1}{S_\Omega}\int_\Omega |T(x,y) - \hat{T}(x, y)|^2 \dd{\Omega}},
\end{equation}
with $S_\Omega$ being area of domain ($S_\Omega = 1$).

The model for Dirichlet problem is
\begin{equation}
    \label{eq:heat-dirichlet-model}
    \hat{T}(x;\theta, \beta) = y\sin\pi x + x(1-x)y(1-y)\cdot N(x;\theta, \beta).
\end{equation}
It is easily verified that the requirements for $G$ and $\tilde{N}$ are 
satisfied. The weight function is actually constructed by the principle 
discussed above: it consists of factors from all four boundary conditions.

\begin{figure}[hptb]
    \centering
    \subfloat[][Analytical solution
        \label{fig:heat-accurate}
        ]{%
	\includegraphics[width=0.45\textwidth]{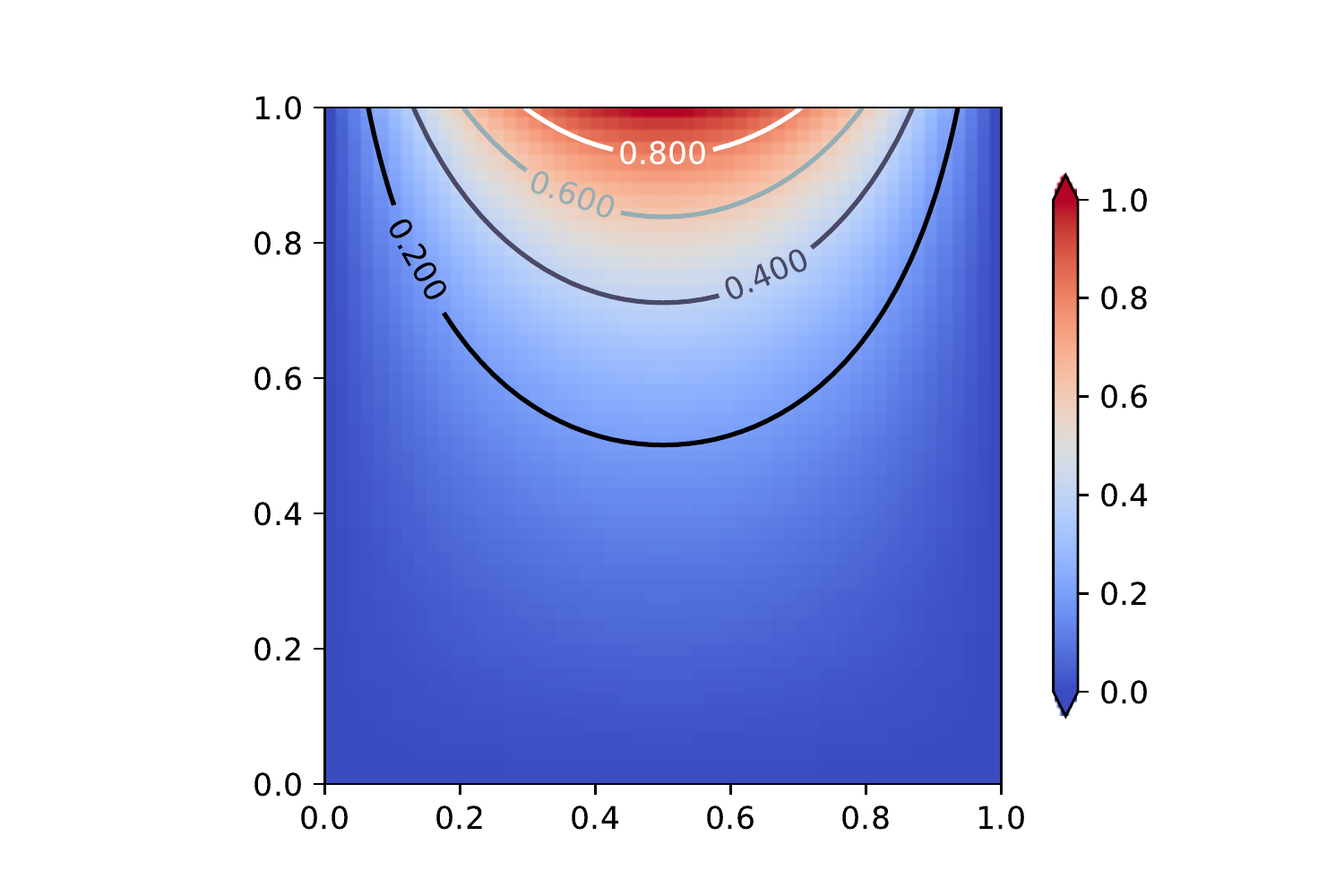}}
    \subfloat[][Numerical Solution
	\label{fig:heat-dirichlet}]{%
	\includegraphics[width=0.45\textwidth]{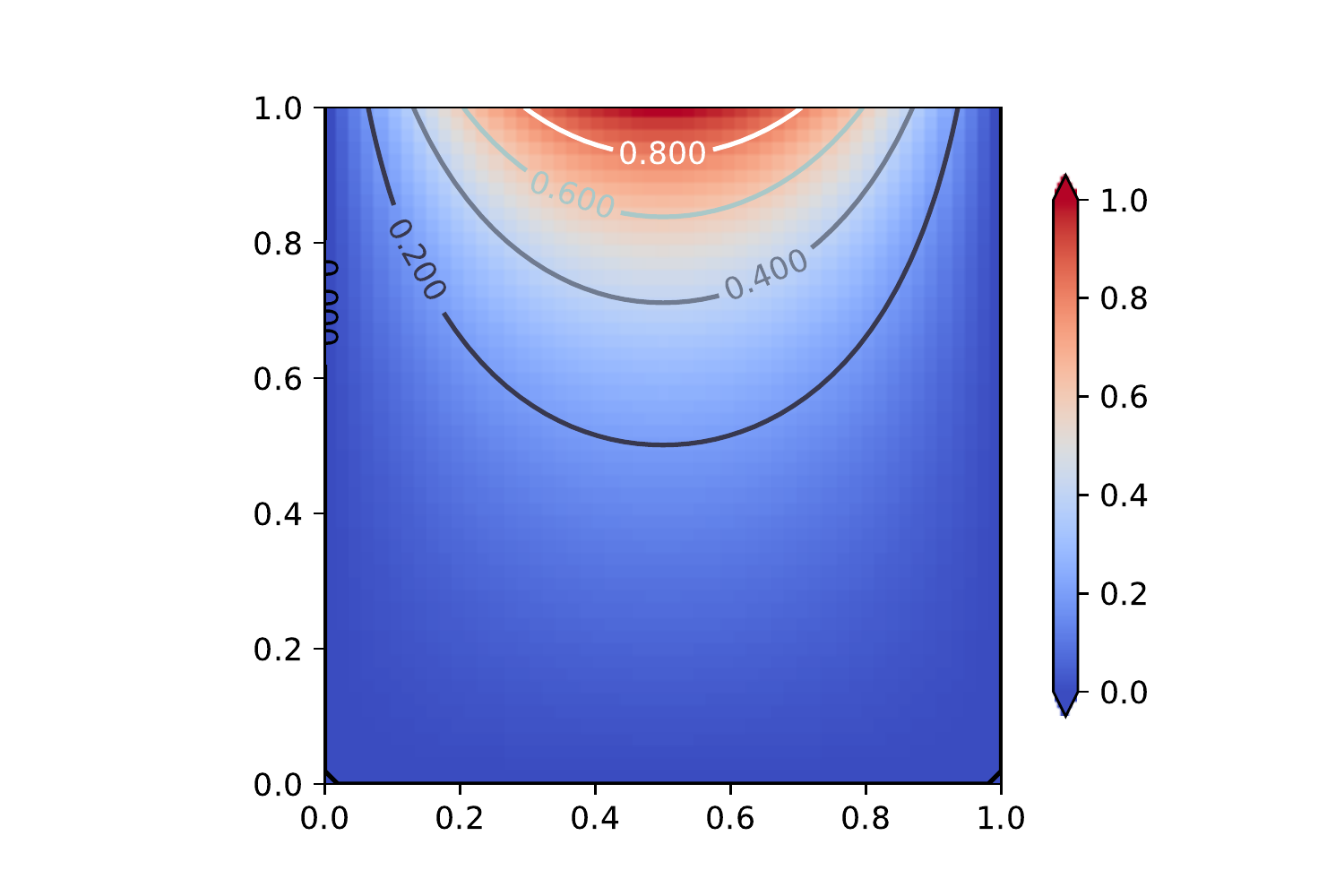}}\\
    \subfloat[][Deviation of bulk term
	\label{fig:heat-error}]{%
	\includegraphics[width=0.45\textwidth]{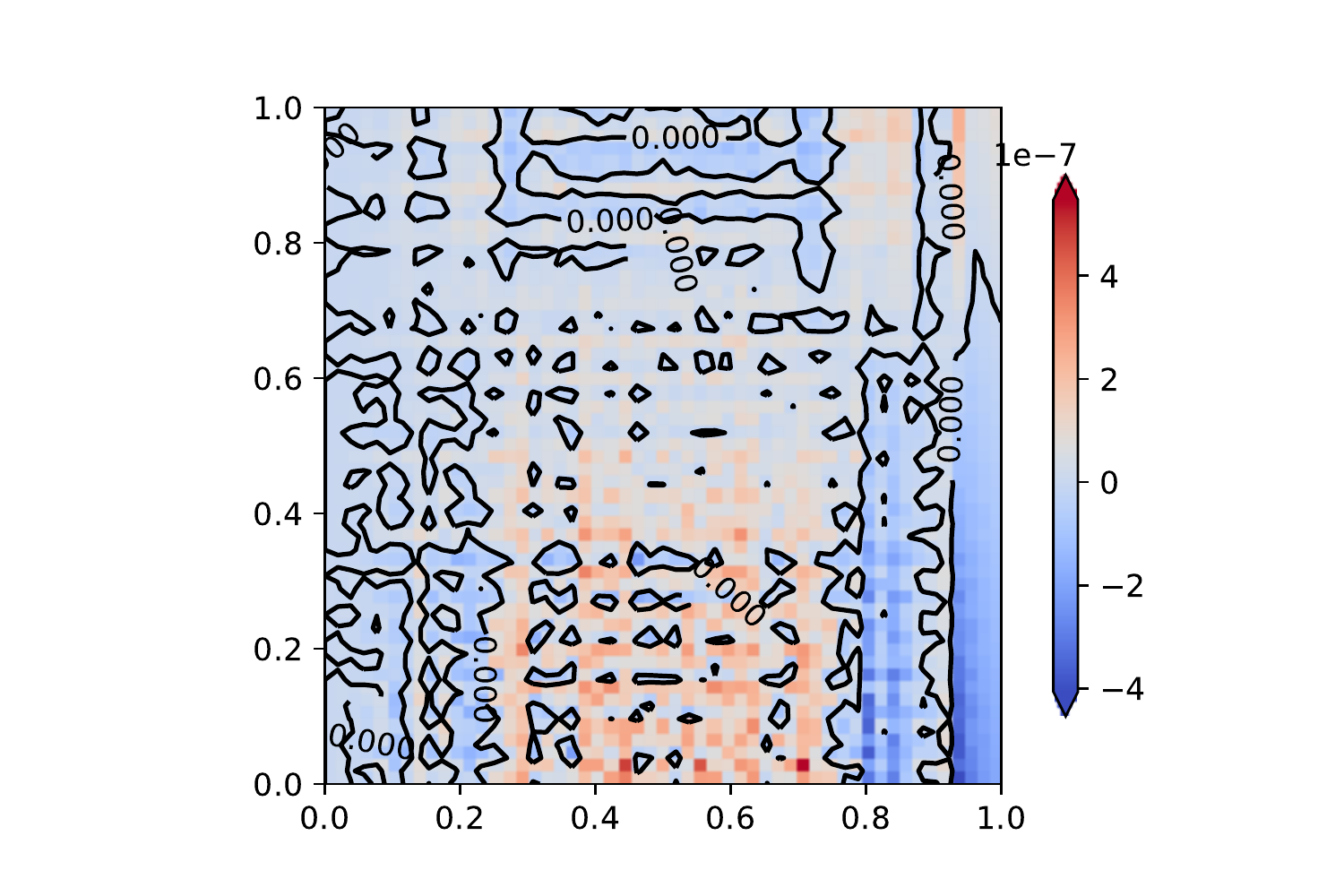}}
    \caption{The results for simulating Laplace equation}
    \label{fig:heat}
\end{figure}

The loss function is constructed similar to \cref{eq:1d-interal-loss}, and the 
simulation is done by a MFN with $2$ hidden layers and $40$ neurons in each 
hidden layer. The data points training the network is a set of $900$ items 
which are vertices of a $30\times30$ uniform mesh on two dimensional unit cube.
The results of simulation is demonstrated in \cref{fig:heat}. 
\Cref{fig:heat-accurate} is contour of analytical solution to 
\cref{eq:heat-solution}, and the simulated solution in 
\cref{fig:heat-dirichlet} matches the analytical solution quite exactly. 
\Cref{fig:heat-error} illustrates the pointwise deviation from analytical 
solution of bulk term. The average error in \cref{fig:heat-dirichlet} is 
$2.8\times10^{-6}$. The similar case calculated by penalty method reported in
\cite{raissi2017physics,liu2019neural} is only of order $10^{-3}$.

\begin{figure}
    \centering
    \subfloat[][Numerical solution
	\label{fig:heat-cd-num}]{%
	\includegraphics[width=0.45\textwidth]{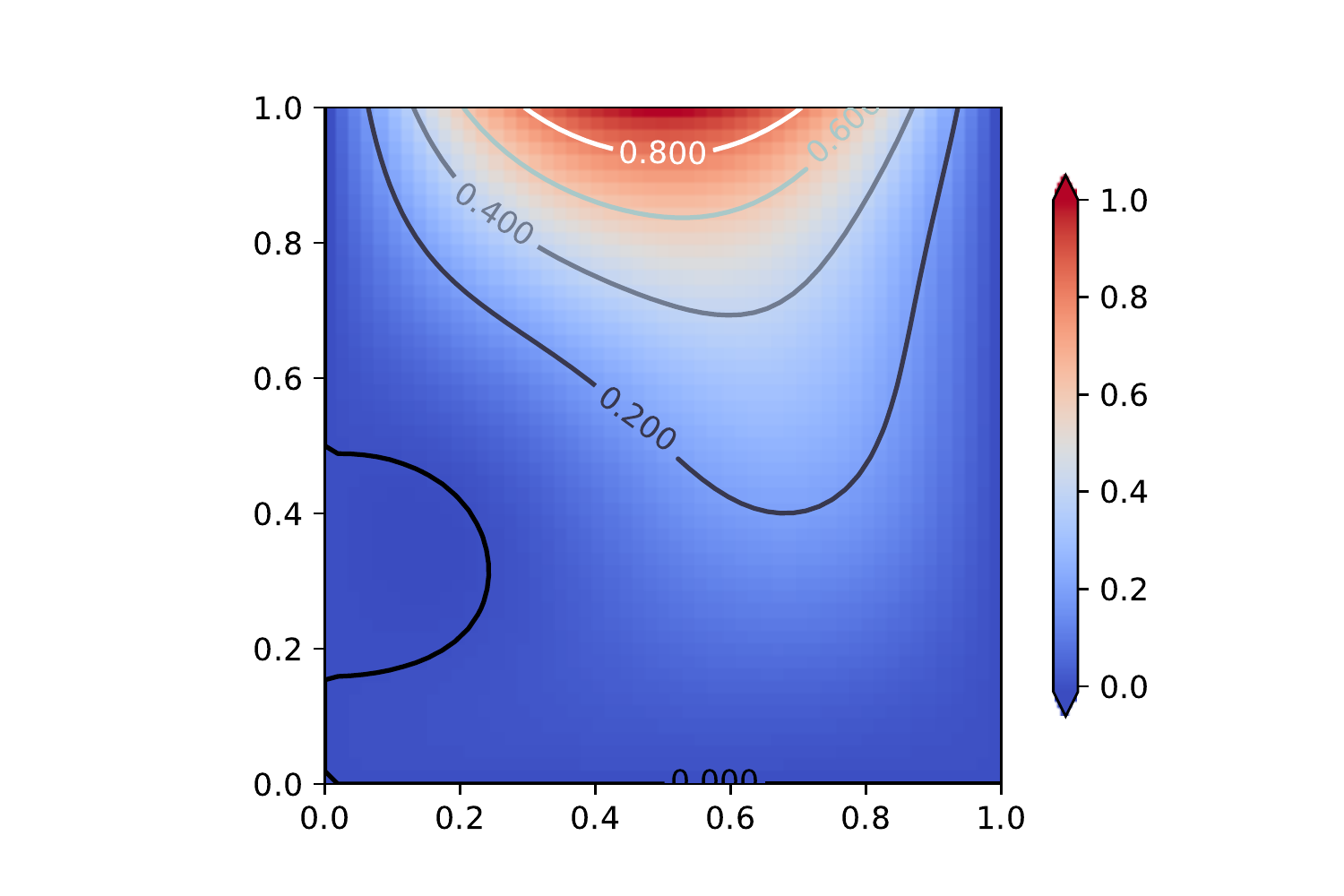}}
    \subfloat[][Analytical solution
	\label{fig:heat-cd-exact}]{%
	\includegraphics[width=0.45\textwidth]{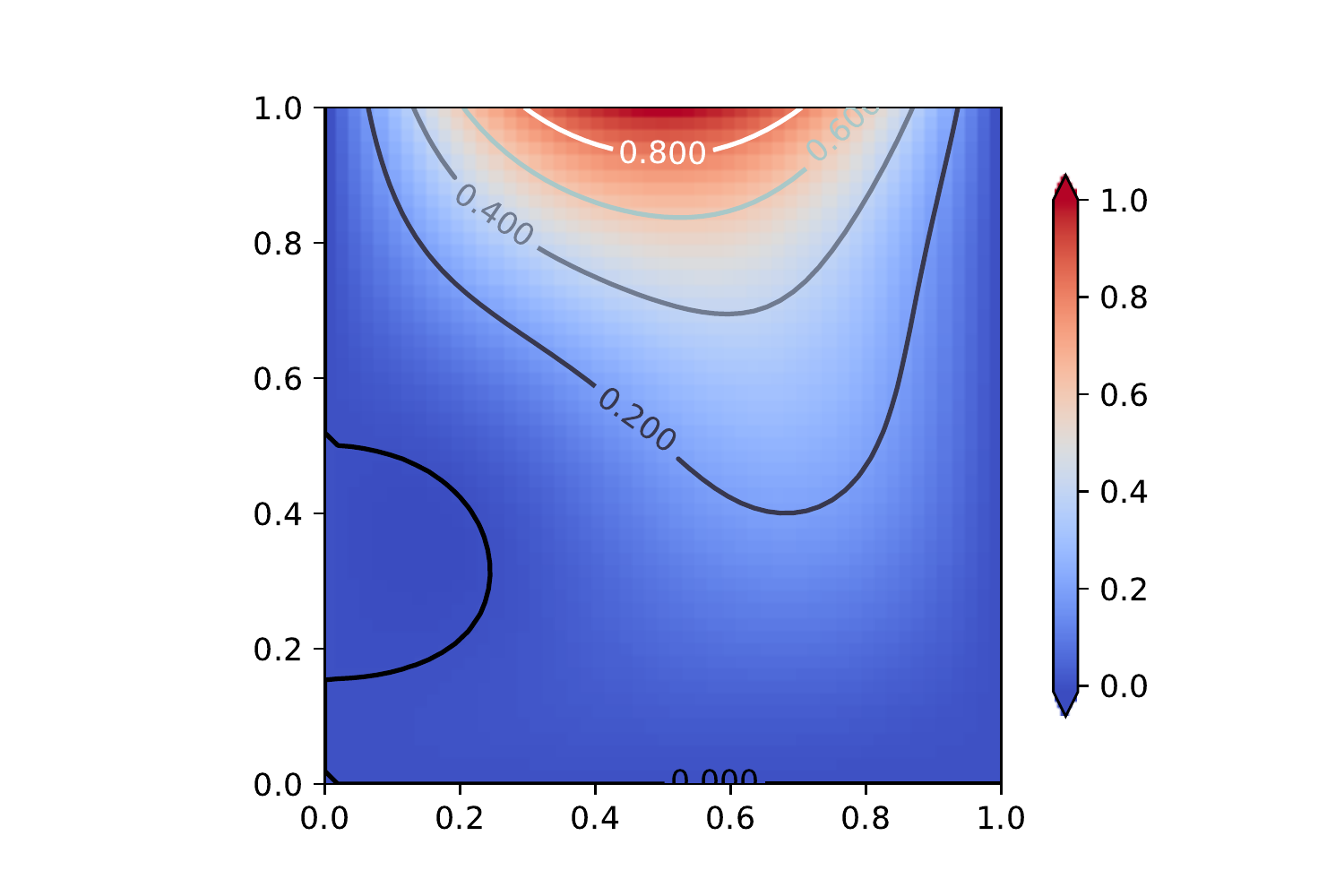}}
    \caption{Convection-diffusion problem with average error of 
	$8.57\times10^{-4}$}
    \label{fig:heat-cd}
\end{figure}

The two properties distinguish CMFN from other methods are its generality and 
accuracy. As for the generality side, this framework is not sensitive to the 
property of differential equation; it works similarly on elliptic, hyperbolic, 
and parabolic problems. One interesting verification is that turning the 
original problem into a convection-diffusion problem:
\begin{equation}
    \label{eq:c-d}
    u\frac{\partial T}{\partial x}+v\frac{\partial T}{\partial y}=\nabla^2 T+f,
\end{equation}
the boundary condition is set the same as in \cref{eq:dirichlet-bd}, as a 
Dirichlet problem, and all other configurations such as network topology and 
training set are all the same as previous. The convection velocities $(u, v)$ 
and source term $f$ are assigned artificially as:
\begin{gather*}
    u(x, y) = y^2 \cos x, \\
    v(x, y) = \frac13 y^3  \sin x, \\
    f(x,y) = y^2\cos x \frac{\pi  \cos \pi  x \sinh \pi  y}{\sinh \pi } + \frac13 y^3\sin x \frac{\pi  \sin \pi x \cosh \pi  y}{\sinh \pi  } + \\ \alpha\left(2\pi^2\sin2\pi x\cos2\pi y-4\pi^2\sin2\pi x\sin^2\pi y\right) -\\
    \alpha \left(2\pi y^2\cos x \cos2\pi x \sin^2\pi y + \frac\pi3 y^3 \sin x \sin2\pi x \sin2\pi y \right)\\
    \alpha = 0.1.
\end{gather*}
This setup ensures that the analytical solution is:
\begin{equation}
    \label{eq:cd-solution}
    T(x,y) =\frac{\sin\pi x\sinh\pi y}{\sinh\pi} - \alpha \sin 2\pi x \sin^2 \pi y
\end{equation} 
The simulated solution is shown in \cref{fig:heat-cd}, its average error is 
$8.57\times10^{-4}$, while the elliptic counterpart of the problem has average 
error of $2.6\times10^{-5}$.

\begin{figure}
    \centering
    \subfloat[][Numerical solution
	\label{fig:heat-pinn-solution}]{%
	\includegraphics[width=0.45\textwidth]{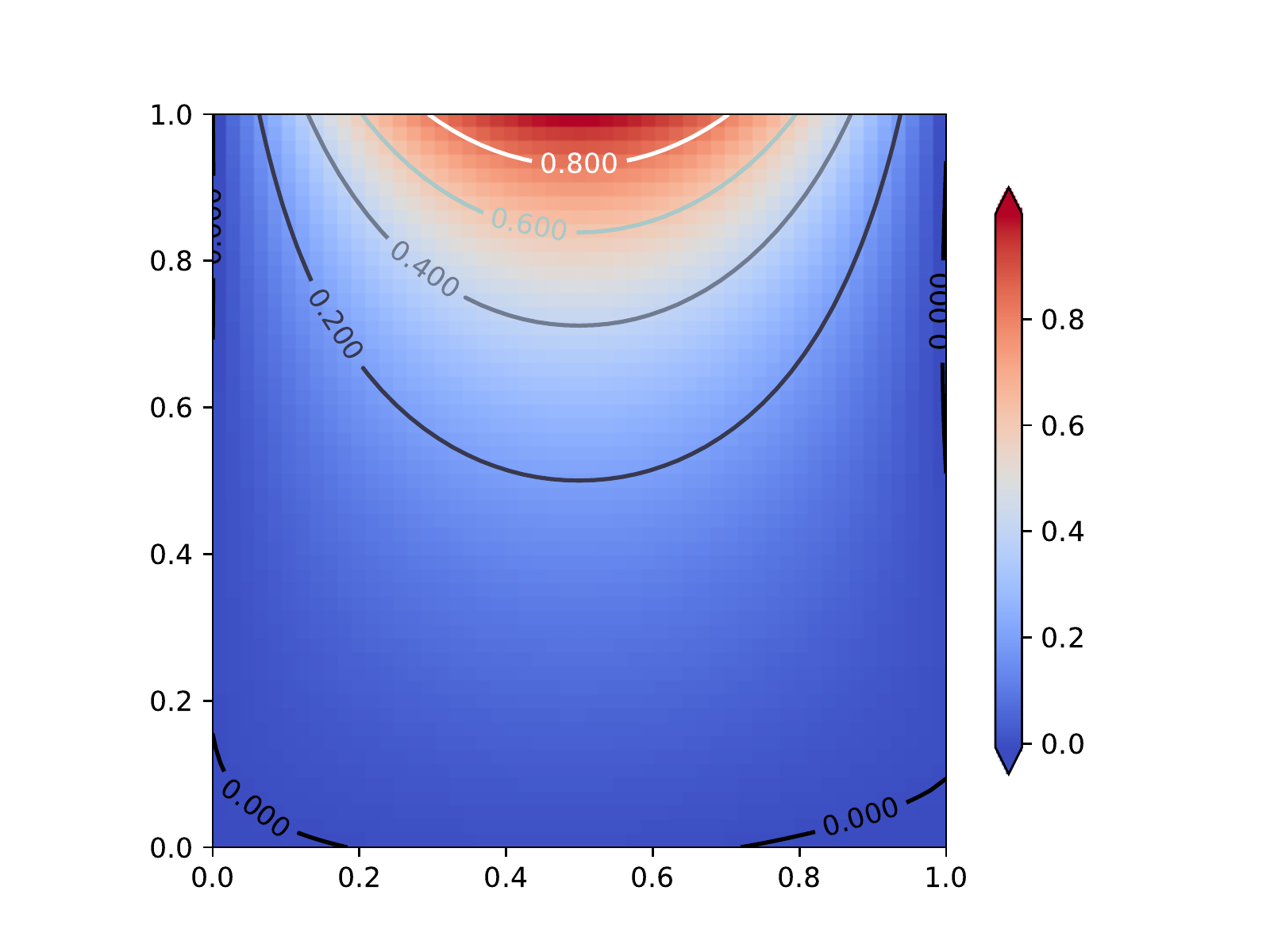}}
    \subfloat[][Error distribution
	\label{fig:heat-pinn-deviation}]{%
	\includegraphics[width=0.45\textwidth]{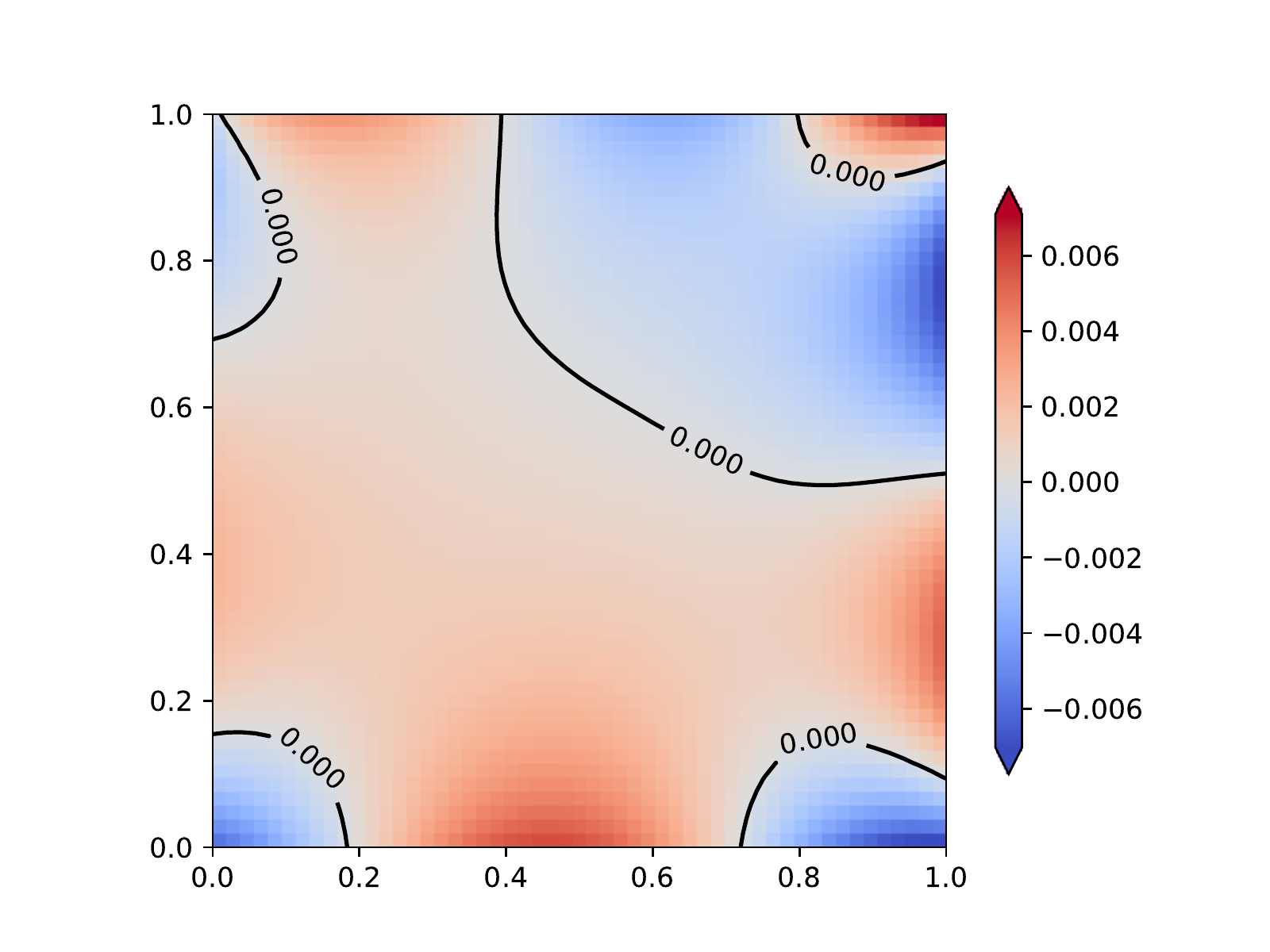}}
    \caption{Numerical solution by PINN method with average error of 
	$1.6\times10^{-3}$}
    \label{fig:heat-pinn}
\end{figure}

Another important property of CMFN is its accuracy. While solving Laplace 
equation, penalty method such as PINN \cite{liu2019neural} is only able to 
achieve an average error of $1.6\times10^{-3}$, and it is easily observed that 
boundary conditions are not accurately satisfied (especially in the two lower 
corners of \cref{fig:heat-pinn}), but CMFN keeps the boundary conditions being 
satisfied accurately in an intrinsic way.

\section{Conclusion and Future Work}

In this paper, we present a novel framework of constructing ODE/PDE solver
based on CMFN method. The numerical method and its application are discussed 
with regard to ODEs and PDEs with Dirichlet boundary conditions.

CMFN method stands out for its generality and accuracy. Traditional neural network
methods based on RBF \cite{mai2001numerical} or penalty methods 
\cite{raissi2017physics} have very limited accuracy. By constructing the trial 
function with a weighted reduced solution as bulk term and a pre-defined 
boundary term, the model satisfies the boundary condition automatically, and as 
a result, the training based on residuals of differential equations could be 
more effective. Moreover, the CMFN method trains the neural network with all 
input data simultaneously, so it is intrinsically able to remain accuracy in 
numerically solving the differential equation on large domain and large time 
span. The iteration methods such as FVM and FDM accumulate truncation error in 
each step, so the scheme has to be designed carefully to be applied to larger
domain and larger time span, while methods based on neural network does not has 
the obsession.

The generality of CMFN framework is in several aspects. The iteration methods 
such as FDM, FVM, and FEM are sensitive to property of PDE, since the growth of 
numerical error differs in hyperbolic, parabolic, and elliptic problems. CMFN 
instead provides a unified method. Compared with traditional neural network
methods based on RBF~\cite{mai2001numerical}, CMFN can be applied similarly on 
both linear and nonlinear problems while the later usually only works on linear 
problem. In this work, very simple network topology (four-layer feedforward 
network with twenty neurons) and very small data (less than $10^3$ points) are 
used, but heat transfer equation and convection-diffusion equation with 
Dirichlet boundary condition on unit cube are solved successfully on the same 
model.

Another property of CMFN method being worth mentioning is the indeterminacy in 
numerical result. In the training stage of our new framework, there is an 
initial guess on MFN\@. Since the network has tremendous parameters, it has to be 
randomly initialized; after training, the loss function would be reduced to a 
small number, but usually not zero, so the optimal solution is not usually 
obtained. The above two factors lead to the result that each specific parameter 
of MFN has rather random behavior. However the overall behavior of the 
computational machine is controllable, because as long as the object function
has enough continuity, the error would reduce along with reduction of loss 
function.

This work is a new starting point in the field of constructing PDE solver for
the authors. There are several works could be considered in the future:
\begin{enumerate}
    \item finding a general method to construct proper form of bulk term for 
	Neumann boundary condition;
    \item finding a systematical method of constructing weight and boundary 
	term, especially for complex geometry;
    \item building larger and deeper networks for more complex problems such as 
	Navier-Stokes equation; and
    \item giving out a more mathematically rigorous proof on existence and 
	uniqueness of reduced solution.
\end{enumerate}

\bibliography{paper.bib}
\bibliographystyle{unsrt}
\end{document}